\documentclass[11pt,a4paper,leqno]{article}
\usepackage[german,english]{babel}
\usepackage[leqno]{amsmath}
\usepackage{amssymb,amsxtra}
\input{diagrams.tex}
\usepackage{latexsym}
\newcommand{\Z}{{\mathbb Z}}       
\renewcommand{\O}{{\cal O}}
\newcommand{\F}{{\mathbb F}}       

\newcommand{\bX}{{\bf X}}
\newcommand{\bY}{{\bf Y}}
\newcommand{\bZ}{{\bf Z}}

\newcommand{\M}{{\cal M}}  
\newcommand{\C}{{\cal C}} 
\newcommand{\cE}{{\cal E}}  
  
\renewcommand{\S}{{\cal S}}  

\renewcommand{\L}{{\cal L}}
\newcommand{\G}{{\tilde{G}}} 
\newcommand{\bG}{{\bf G}}

\newcommand{\tG}{{\tilde{G}}}
\newcommand{\PR}{{\cal PR}}
\renewcommand{\H}{{\cal H}}
\newcommand{\mlabel}[1]{\label{#1}}

\newcommand{\brho}{{\boldsymbol{\rho}}}
\newcommand{\bbe}{{\boldsymbol{\beta}}}

\newcommand{\arrsim}{\stackrel{\sim}{\longrightarrow}}

\newcommand{\Res}{\mbox{{\rm Res}}}
\newcommand{\coker}{\mbox{{\rm coker}}}
\newcommand{\Ind}{\mbox{{\rm Ind}}}
\newcommand{\Inf}{\mbox{{\rm Inf}}}
\newcommand{\Fix}{\mbox{{\rm Fix}}}
\newcommand{\Hom}{\mbox{{\rm Hom}}}
\newcommand{\End}{\mbox{{\rm End}}}
\newcommand{\Jac}{\mbox{{\rm Jac}}}
\newcommand{\soc}{\mbox{{\rm soc}}}
\newcommand{\stab}{\mbox{{\rm stab}}}
\newcommand{\im}{\mbox{{\rm im\,}}}
\newcommand{\hd}{\mbox{{\rm hd}}}

\newcommand{\R}{\mbox{{\rm R}}}
\newcommand{\T}{\mbox{{\rm T}}}
\renewcommand{\mod}{{\mathfrak{Mod}\,}}
\renewcommand{\lat}{{\mathfrak{Lat}\,}}
\newcommand{\bker}{{\bf ker}}

\newcommand{\be}{\beta}

\newcommand{\D}{{\cal D}}
\newcommand{\tD}{\tilde{\cal D}}
\newcommand{\E}{{\cal E}}
\newcommand{\Or}{{\cal O}}
\newcommand{\cN}{{\cal N}}

\newcommand{\cB}{{\cal B}}

\newcommand{\cD}{{\cal D}}

\newcommand{\cX}{\mathfrak{X}}

\newcommand{\cL}{{\cal L}}

\newcommand{\NLM}{N_{L\subseteq M}}

\newcommand{\NLGd}{N_{L\subseteq G'}}
\newcommand{\NMGd}{N_{M\subseteq G'}}
\newcommand{\tilNLG}{\tilde{N}_{L\subseteq G}}
\newcommand{\tilNLM}{\tilde{N}_{L\subseteq M}}
\newcommand{\cNNML}{\cN_{M,L}^G}
\newcommand{\cNNGdL}{\cN_{G',L}^G}
\newcommand{\cnorm}[3]{\cN_{#1,#2}^{#3}}

\newcommand{\Rs}{R\in\{K,\Or,F\}}

\newcommand{\pf}{\smallskip {\bf Proof\hspace{1em}}}
\renewcommand{\simeq}{\cong}

\def\epf{\ifmmode\eqno\Box\medskip\else{\unskip\nobreak\hfil%
  \penalty50\hskip2em\hbox{}\nobreak\hfil$\Box$
  \parfillskip=0pt \finalhyphendemerits=0\penalty-100\medskip}\fi}

\makeatletter
\@addtoreset{equation}{section}
\def\theequation{\thesection.\arabic{equation}}
\makeatother
\newcommand{\rmref}[1]{{\rm\ref{#1}}\it}

\newcommand{\fC}{\mathfrak C}

\newtheorem{thm}[equation]{Theorem}
\newtheorem{lemma}[equation]{Lemma}
\newtheorem{notation}[equation]{Notation}

\newtheorem{cor}[equation]{Corollary}

\newtheorem{rem}[equation]{Remark}

\newtheorem{hypo}[equation]{Hypothesis}

\newtheorem{summ}[equation]{}
\newtheorem{ddef}[equation]{Definition}

\makeatletter

\def\enumerate{%
  \ifnum \@enumdepth >3 \@toodeep
  \else\advance\@enumdepth \@ne
       \edef\@enumctr{enum\romannumeral\the\@enumdepth}%
       \topsep\z@\parskip\z@\partopsep\z@                                     
       \list
       {\csname label\@enumctr\endcsname}
       {\@nmbrlisttrue\setcounter{\@enumctr}{0}\let\@listctr\@enumctr
       \def\makelabel##1{\hss\llap{##1}}
       \parsep\z@\itemsep\z@}%
  \fi}

\def\Number#1{\refstepcounter{equation}%
              \leqno(\theequation)\if*#1\else\label{#1}\fi}     

\makeatother

\topmargin--5mm
\begin{document}
\thispagestyle{empty}

\title{Generalized $q$-Schur algebras and modular representation\\
  theory of finite groups with split $(BN)$-pairs}
\author{
            Richard Dipper and Jochen Gruber\\
            Mathematische Institut B\\
            Universit{\"a}t Stuttgart\\
            Postfach 80 11 40\\ 
            70550 Stuttgart \\
            Deutschland
}
\makeatletter
\let\@makefnmark\relax
\footnotetext{\kern-6mm%
\begin{tabular}{ll}
E--mail:& rdipper@mathematik.uni--stuttgart.de\\ 
\              & jgruber@mathematik.uni-stuttgart.de
\end{tabular}\\
{\it 
A.M.S. subject classification (1991): 20C20, 20C33, 20G05, 20G40}\\
This paper is a contribution to the DFG project on ``Algorithmic
number theory and algebra''. The authors acknowledge support
from the Isaac Newton Institute, Cambridge}
\makeatother

\maketitle 
 
\begin{abstract} We introduce a generalized version of a 
$q$-Schur algebra (of parabolic type) for arbitrary Hecke algebras
over extended Weyl groups. We
describe how the decomposition matrix of a finite group with split $BN$-pair, 
with respect to a non-describing prime, can be
partially described by the decomposition matrices of suitably chosen 
$q$-Schur algebras. We show that the investigated structures occur naturally
in finite groups of Lie type. \\

\end{abstract}

{\bf\Large Introduction}\mlabel{intro}
\\

\noindent In a series of papers (\cite{di2,di3,ja,dija1} Gordon James and the 
first named author used representations of Hecke- and of $q$-Schur algebras to 
derive a classification of the $\ell$-modular irreducible representations 
of general linear groups $G=GL_n(q)$ for primes $\ell$ which are coprime 
to $q$. It turned out that the $\ell$-decomposition matrices  of $G$ are 
completely determined by decomposition matrices of certain $q$-Schur algebras
by an algorithm involving only combinatorics as partitions 
and the Littlewood-Richardson rule (see \cite{dija1}). 

The connection between the group algebra of $G$ on the one side, and 
Hecke- and $q$-Schur algebras on the other side is given by certain 
functors $H$, called quotients of $\Hom$-functors and their preinverses, which 
were investigated in \cite{di,di4}. First Hecke algebras come up as 
endomorphism rings of certain induced modules for $G$, and these functors 
connect their representations with those of $G$. This part of the theory 
was extented to arbitrary finite groups of Lie type in \cite{difl1,difl2}.
For general linear groups $q$-Schur algebras turn up as endomorphism rings
of certain Hecke algebra modules. In fact quotients of Hom-functors can be 
defined here as well, one obtains a $q$-analogue of the classical Schur 
functors. The functors $H$ are then used to show that the $q$-Schur algebras
are isomorphic to endomorphism rings of certain $G$-modules as well, and 
that one can define again quotients of Hom-functors $S$, now connecting 
representations of $q$-Schur algebras and of $G$. In \cite{grhi} Hiss and the 
second named author extended these results to finite classical groups 
in the special case of linear primes $\ell$ (about a third of the primes
dividing the group order).  

The main purpose of this paper is to exhibit how far these methods carry in 
the general case of finite groups $G$ of Lie type for arbitrary primes 
$\ell$ different from the describing characteristic of the group. 
One of our main results is the extension of the second step, the 
$q$-Schur algebra approach in the situation of \cite{difl1,difl2}. In 
addition we shall show that the methods work in other cases too. Here 
a fundamental step has been done by Geck, Hiss and Malle, who showed
that the endomorphism ring of Harish-Chandra induced cuspidal irreducible 
modules is always a Hecke algebra (see \cite{gehima2}), extending 
Howlett-Lehrer theory \cite{hole1}. In addition they constructed corresponding
quotients of Hom-functors.

However it should be pointed out that our results are not as nice and complete 
as in the case \cite{dija1} of general linear or \cite{grhi} of classical 
groups for linear primes: In general we can construct only a part of the
irreducible representations of $G$ and we get only partial information on 
the decomposition numbers in terms of generalized $q$-Schur algebras. 
It seems to be also clear, that the information provided by our results 
here is about everything what $q$-Schur algebras can give. To obtain 
information on the remaining irreducible representations of $G$ and the 
corresponding decomposition numbers one probably needs fundamentally new 
methods. 

We now describe briefly how we proceed in this paper. 
If $G$ is a finite group with split $BN$-pair, one has the
distribution of complex characters in so called Harish-Chandra series,
HC-series for short. They are indexed by pairs whose first entry is a 
representative $L$ of conjugacy classes of Levi subgroups of $G$, and the 
second entry is a cuspidal irreducible character of $L$. Attached 
to such an HC-series we have a Hecke algebra and one can define functors 
between the categories of $G$-modules  and the module
categories of the respective Hecke algebras. 

In section one we first review the results on Hom-functors needed in
the following. Then we define quotients of projective $\O G$-lattices 
by factoring out constituents which do not belong to a fixed HC-series. 
These quotients are defined by a functor which is compatible with 
HC-induction.

In section two we introduce the central object of the paper, 
projective restriction systems associated with HC-series. It turns out that 
this axiomatic setting enables us to step up from using 
Hecke- to using $q$-Schur algebras to obtain  
information on representations belonging to the corresponding  
HC-series. We investigate how modules, which are contained in such a
projective restriction system, behave under Harisch-Chandra restriction. 
We exhibit in particular a close relation between the poset
structure of Levi subgroups of $G$ and module structure of the modules
in the projective restriction system. As a special case one can bear
in mind the Steinberg lattices of finite groups of Lie type.

In section three we then connect our results of section two and those
of section one to get our main result on decomposition numbers.

In section four we consider the most important groups with split
$BN$-pairs, the finite groups of Lie type. We show that projective
restriction systems appear regularly for these groups 
and apply our results to obtain at the end a new concrete result for 
the unipotent characters of finite unitary groups.

This paper has a long history: A preliminary version was out already at the 
end of 1996 and it was referred to in several papers which meanwhile appeared. 
In between there were several improvements and generalisations, which may 
explain the delay in part, (however the discovery of several gaps in proofs, 
which had to be closed, contributed subtantially as well). The main results 
were announced at several conferences, in \cite{dghm1} and in \cite{gr3}.

We wish to thank the organizers of the half year program on 
algebraic groups and related finite groups at the Isaac Newton Institute 
in Cambridge in the first half of 1997 and the University of Illinois at
Chicago. During an extended visit at the Isaac Newton Institute by the first 
and an one year visit in Chicago in 1996/97 by the second author most of this 
article was written up.   

\section{Preliminaries}\mlabel{prel}
Throughout $G$ denotes a finite group. We let
$\O$ be a complete discrete valuation ring. The quotient field of $\O$
is denoted by $K$ and its residue field by $k$. We assume that $k$ is
of characteristic $\ell$ for some prime $\ell$ dividing the order of
$G$. Moreover we take $(K,\O,k)$ to be an $\ell$-modular
splitting system for $G$. Thus both fields, $K$ and $k$ are splitting
fields for all subgroups of $G$.

For $H\leq G$ we denote the normalizer of $H$ in $G$ by $\cN_G(H)$.

Modules are, if not stated otherwise, right modules, and homomorphisms
act from the opposite side. Let $R$ be a commutative domain. An
\emph{$R$-order} $A$ is an $R$-algebra which is free and finitely
generated as $R$-module, and an \emph{$A$-lattice} is an
$A$-module which is finitely generated and free as $R$-module. The
category of $A$-lattices is denoted by $\lat_A$ and of $A$-modules by
$\mod_A$. Since $R$ is
commutative we may write the action of scalars on $A$-modules on the
left assuming always that $R$ acts centrally on every module.

If possible we omit tensor symbols. For instance, if $B$ is some 
sublattice of $A$ and $M$ is a $B$-lattice we frequently write 
$MA$ for the induced module $M\otimes_B A$. Or if $R=\O$, and if
$M$ is an $A$-lattice, $KA$ denotes the finite dimensional
$K$-algebra $K\otimes_\O A$, and $KM$ the $KA$-module $K\otimes_\O M$.
The $k$-Algebra $k\otimes_\O A$ is denoted by $kA$ or by $\bar{A}$
and the $\bar{A}$-module $k\otimes_\O M$ by $\bar{M}$.
Similarly, if $H$ is a group and $\chi$ is a $KH$-character,
its associated Brauer character, that is the restriction of $\chi$ to
$\ell$-regular classes of $H$, is denoted by $\bar{\chi}$. 

The \emph{radical} of an $A$-lattice $M$ is the intersection of its
maximal sublattices and is denoted by $\Jac(M)$. The \emph{socle}
$\soc(M)$ of $M$ is the maximal completely reducible submodule of $M$
and the \emph{head} $\hd(M)$ of $M$ is the maximal completely reducible
factor module of $M$. Note that $\hd(M) = M/\Jac(M)$.


We need some results on quotients of Hom-functors from \cite{di} and
\cite{di4}. We begin with the basic setup there and supplement a few
further results which will be needed later on. 

Let $T$ be a semiperfect $R$-algebra which is finitely
generated as $R$-module. All occurring modules are, if not stated
otherwise, finitely generated. Let $M \in \mod_T$, and let $\beta : P \to
M$ be a projective presentation of $M$. Thus $P \in \mod_T$ is
projective, and $\beta$ is an epimorphism. Let $\E = \End_T(P)$.
We take $\E_{\beta} = \{ \phi \in \E\, |\, \phi (\ker\beta )\subseteq 
\ker\beta\}$.
The endomorphism ring of $M$ is denoted by $\H$. Obviously
$J_{\be} = \{\psi\in\E\, |\,\im\psi \leq \ker\beta\}$ is an ideal of $
\E_{\be}$ and $\E_{\be}/J_{\be} \cong \H$ as $R$-algebra canonically, 
(comp. \cite[2.1]{di}). Our basic hypothesis is now the following:

\begin{hypo}\label{hypo}
For $M \in \mod_T$ we say that the projective presentation $\be : P
\rightarrow M$ of $M$ satisfies \rmref{hypo}, if $\E_{\be} = \E$.
\end{hypo}

Then $J_{\be}$ is an ideal of $\E$, and we identify $\H$ and
$\E /J_{\be}$ by the canonical isomorphism induced by $\be$. We have
now a functor 
$$
H = H^\be = H^\be_M: \mod_T\rightarrow \mod_{\H}
$$ 
which takes the
$T$-module $V$ to the $\H$-module $\Hom _T(P,V)/\Hom _T(P,V)J_{\be}$.
On maps $H^\be$ is defined in the obvious way.

\begin{notation}\label{tauandt} {\rm The \emph{$P$-torsion submodule}
$t_P(V)$ is the unique maximal submodule $X$ of $V \in \mod _T$ with
respect to the property that $\Hom _T(P,X) = (0)$. The \emph{kernel}
$\bker P$ of $P$ is the full subcategory of $\mod _T$, whose objects are
$T$-modules $V$ with $\Hom _T(P,V) = (0)$. So $V \in \bker P$ precisely
if $t_P(V) = V$. We have a functor $A_P: \mod _T \to \mod _T$ taking
$V \in \mod _T$ to $V/t_P(V)$. We say that $V$ is \emph{$P$-torsionless}
if $t_P(V) = (0)$. So $A_P(V)$ is the maximal $P$-torsionless factor
module of $V$. We define the functor $\hat H_M = \hat H: \mod _{\H}
\to \mod _T$ to be the composite functor $A_P \circ (-\otimes _{\H}
M)$. For $V \in \mod _T$ the \emph{trace} or \emph{combined image}
$\tau_P(V)$ of $P$ in $V$ is
the $T$-submodule of $V$ generated by the images of all homomorphisms
$\phi : P \to V$.}
\end{notation}

The following lemma follows immediately from projectivity of $P$:

\begin{lemma}\label{tauundtfunct} Let $P,V,U\in\mod_T$ and suppose
  that $P$ is projective. Let $f:V\to U$ be $T$-linear. Then 
$$
f(\tau_P(V))\leq \tau_P(U)\quad\text{and}\quad f(t_P(V))\leq t_P(U).
$$
\end{lemma}

Thus $\tau_P$ and $t_P$ are really endofunctors of $\mod_T$.
 
It was shown in \cite[2.16]{di} that $\hat H$ is a right inverse of
the functor $H$.  Thus in particular, if $T$ is a finite dimensional
algebra over some field, $\mathfrak H$ is a complete set of non-isomorphic
irreducible $\H$-modules, and $\mathfrak T_M$ a complete set of
non-isomorphic irreducible $T$-modules which are not taken to the zero
module under the functor $H$, then  $H$ induces a bijection
between $\mathfrak T_M$ and $\mathfrak H$ \cite[2.28]{di}.

More precisely we have the following theorem (see \cite[1.2]{di4}):

\begin{thm}\label{1.3}
Suppose that $\be:P\to M$ satisfies \rmref{hypo}. Then 
\begin{enumerate}
\item If $M$ is $P$-torsionless and $X$ is a right ideal of $\H$. Then 
$$
\hat{H}(X) = XM \cong A_P(X\otimes_\H M).
$$
\item If $R$ is a field then $\mathfrak{T}_M$ is a complete set of
  nonisomorphic irreducible constituents of $\hd(M)$. Then every
  indecomposable direct summand of $M$ has simple head and factoring out
  the Jacobson radical induces a bijection between the isomorphism
  classes of indecomposable direct summands of $M$ and the elements of
  $\mathfrak{T}_M$. 

\end{enumerate}\end{thm}

Let now $T$ be an $\O$-order in the semisimple $K$-algebra $KT$. We
say $V$ is an \emph{irreducible $T$-lattice} if $KV = K\otimes_\Or V$
is an irreducible $KT$-module.

Recall that an $kT$-module $V$ is \emph{liftable}, if $V =
k\otimes_\Or U$ for some $T$-lattice $U$. Moreover, if 
$U, V$ are $T$-lattices,  then
\begin{equation}
\Hom_{KT}(KU,KV) \cong K\otimes_{\Or}\Hom_{T}(U,V)
\end{equation}
canonically. Similarly  we have an homomorphism of $k$-spaces
\begin{equation}\label{1.6.2}
\overline{\Hom_T(U,V)}=k\otimes_{\Or}\Hom_T(U,V) \to 
\Hom_{\bar{T}}(\bar{U},\bar{V})
\end{equation}
which is injective but not necessarily surjective. The homomorphisms 
in the image of that map are called \emph{liftable}, and similarly
$\Hom_{\bar{T}}(\bar{U},\bar{V})$ is said to be \emph{liftable}, if this
homomorphism is bijective.

We take $M\in\lat_T$ and assume that the projective presentation 
$\be : P\to M$ satisfies (1.1). By \cite[4.7]{di} Hypothesis \ref{hypo} 
(for $R = \Or $) is
equivalent to $\ker(1_K\otimes_\O\be)$ having no irreducible
constituent in common with $KM$. Moreover if this holds, the
decomposition matrix of $\H$ is a submatrix of the decomposition
matrix of $T$, \cite[4.10]{di}. In \cite[4.7]{di} it was shown that 
\ref{hypo} holds then for $1_R\otimes_\Or\be$ for all choices of 
$R \in \{K,\Or,k\}$. The three
resulting Hom-functors and their right inverses are distinguished by a
suffix $R$. We summarize:

\begin{thm}\label{decmatrix} Let $\be:P\to M$ be a projective
  presentation of the $T$-lattice $M$ and suppose that $K\ker(\be)$
  and $KM$ have no irreducible constituent in common. Then
  $1_R\otimes_\O\be$ satisfies Hypothesis \rmref{hypo} for all choices
  of $R\in\{K,\O,k\}$, the endomorphism ring of $\bar{M}$ is liftable
  and the decomposition matrix of $\H=\End_T(M)$ is part of the
decomposition matrix of $T$.
\end{thm}

\begin{lemma}\label{ptl} Let $\be:P\to M$ be a projective presentation
of the $T$-lattice $M$ satisfying Hypothesis \rmref{hypo}. Then $RM$
is $RP$-torsionless for $R = K$ and $\Or$.
\end{lemma}

Recall that an algebra $A$ over some field is called \emph{quasi
  Frobenius}, if there exists a nondegenerate, associative
bilinearform on $A$. This implies in particular, that the regular
representation of $A$ is injective, that is $A$ is a self-injective
algebra, and every homomorphism between two right ideals is given by
left multiplication by some element of $A$.

If the form is in addition symmetric, we call $A$
\emph{symmetric}. For instance, group algebras over fields are
symmetric algebras. Here is a generalisation to $\O$-orders:

\begin{notation}\label{isym} {\rm An $\O$-order $\H$ is called \emph{integrally
quasi Frobenius} if there exists an associative bilinearform on $\H$ 
whose Gram determinant with respect to any basis of $\H$ is a unit in
$\O$.

A sublattice $V$ of some $\O$-lattice $M$ is called \emph{pure} if the
factor module $M/V$ is a torsion free (and hence free)
$\O$-module. If $V\leq M$ then the intersection of all pure
sublattices of $M$ which contain $V$ is the unique minimal pure
sublattice containing $V$ and is denoted by $\sqrt{V}$. One checks
easily that $\sqrt{V}=KV\cap M$.}
\end{notation}

Note that the bilinearform on $\H$ induces nondegenerate
associative bilinearforms on $R\H$ for $R=K$ and $R=k$, such that
these algebras are quasi Frobenius. In our applications later on we
shall deal with Hecke algebras which are known to be integrally quasi
Frobenius (indeed integrally symmetric).

Here is a result, which we will need later. For a proof see 
\cite[1.30 and 1.36]{di4}:

\begin{thm}\label{hom1corr} Let $\be:P\to M$ be a projective
  presentation of the $T$-lattice $M$ satisfying Hypothesis
  \rmref{hypo}. Let $X$ and $Y$ be pure right
ideals of $\H= \End_T(M)$. Suppose that $\H$ is integrally quasi
Frobenius. Denote the 
associated Hom-functor by $H=H^\be$. Then \begin{enumerate}
\item $H(XM)=H(\sqrt{XM})=X$.
\item Every homomorphism from $\sqrt{XM}$ to $\sqrt{YM}$ maps $XM$
  into $YM$. Restricting homorphisms induces an isomorphism
$$
\Hom_{RT}(R\otimes_{\Or}\sqrt{XM},R\otimes_{\Or}\sqrt{YM})\arrsim
\Hom_{RT}(XM,YM)
$$
for $R=\O,K$.
\item The functor $H$ induces an isomorphism
\begin{equation}\label{homlatt}
H_R:  \Hom_{RT}(R\otimes_{\Or}\sqrt{XM},R\otimes
_{\Or}\sqrt{YM}) \arrsim\Hom_{R\H}(RX,RY),
\end{equation}
for $R = \Or$ and $K$.
\end{enumerate}  
\end{thm}


The following result is a refinement of \cite[3.2 and 3.3]{difl1}.

\begin{lemma}\label{homdim} Let $G$ be a finite group and let $(K,\O,
k)$ be a split $\ell$-modular system for $G$ for some prime $\ell$. 
Let $P$ be a projective $\O G$-lattice and suppose that the irreducible
$KG$-module $M$ occurs with multiplicity one as constituent of
$KP$. Let $X$ be an $\O G$-lattice in $M$ such that $KX=M$. Then
$$
\Hom _{\O G}(P,X)\simeq \O,
$$
and for a generator $\varphi$ of $\Hom_{\O G}(P,X)$ we have the
following: The projective cover of the image $V=\im(\varphi)$ is an 
indecomposable direct
summand $Q$ of $P$. The simple $kG$-module $D=\bar{Q}/\Jac\bar{Q}$ is
the head of $\bar{V}$, that is $\bar{V}/\Jac\bar{V}\simeq D$, and $D$ occurs
as composition factor of $\bar{V}$ with multiplicity one. Moreover,
$$
\End_{kG}(\bar{V})\simeq k\otimes _{\O}\End_{\O G}(V)\simeq k
$$
and the image $U$ of $\bar{\varphi}=1_k\otimes _{\O}\varphi:\bar{P}
\rightarrow\bar{X}$
is an epimorphic image of $\bar{V}$ and $\bar{Q}$ is the projective
cover of $U$. Thus $U$ has simple head $D$ occurring with multiplicity
one as composition factor of $U$. In particular
$$
\End _{kG}(U)\simeq k.
$$
\end{lemma}
\pf Let $P=P_1\oplus\ldots\oplus P_m$ be a decomposition of
$P$ into a direct sum of projective indecomposable $\O
G$-lattices. Then $KP=KP_1\oplus\ldots\oplus KP_m$ and $KX$ occurs as
irreducible direct summand of precisely one $KP_i$, say of $KP_1$, and
we set $Q=P_1$. Thus restricting maps induces an isomorphism 
$$
\Hom _{KG}(KP,KX)\simeq \Hom _{KG}(KQ,KX),
$$
where the left hand side is one dimensional. Since 
$$
\Hom _{KG}(KP,KX)\simeq K\otimes _{\O}\Hom _{\O G}(P,X)\simeq\Hom _{KG}(KQ,KX)
$$
we conclude that 
$$
\Hom _{\O G}(P,X)\simeq \O\simeq \Hom_{\O G}(Q,X).
$$
The generator $\varphi$ of $\Hom _{\O G}(P,X)$ is a generator of $\Hom
_{\O G}(Q,X)$ too, hence $\varphi:Q\to\im(\varphi)=V$ is the minimal
projective cover of $V$. Being indecomposable $\bar{Q}$ has simple
head $\hd(\bar{Q}) = \bar{Q}/\Jac(\bar{Q})$ and this must be the head 
of $\bar{V}$ as well. By Brauer's reciprocity law the
multiplicity of $D$ as composition factor of $\bar{X}$ equals the
multiplicity of $M=KX$ as irreducible constituent of $KQ = KP_1$, and
hence is one. Thus $D$ occurs only once as composition factor of
$\bar{V}$ as well. If $\wp$ denotes the generator of the unique
maximal ideal $(\wp)$ of $\O$, then $k=\O/(\wp)$ and we have
$\bar{V}\cong V/\wp V$, whereas the image $U$ of $\bar{\varphi} =
1_k\otimes_\O\varphi:\bar{P}\to\bar{X}$ is given as  $U = V/(V\cap\wp
X)$. Since $\wp V\subseteq \wp X$, the module $U$ is an epimorphic
image of $\bar{V}$ and the remaining claims of the lemma follow.
\epf

\begin{lemma}\label{inprel2} Let $G$ and $(K,\O,k)$ be as in the
  previous lemma. Let $X$ be an $\O G$-lattice and suppose that its
  character $\chi$ decomposes $\chi = \chi_1+\chi_2$, where
  $\chi_2$ is the character of a pure sublattice $M$ of $X$ and $\chi_1$
  is an irreducible character afforded by $S = X/M$. Assume further
  that $KS$ occurs in $KP$ with multiplicity one, where $P$ is the
  projective cover of $S$ and let $D = \bar{S}/\Jac(\bar{S})$. Then
  $D$ is an irreducible $kG$-module and occurs with multiplicity one in
  $\bar{S}$. Moreover
$$
\End_{kG}(\bar{S}) \cong k
$$
and, if $D$ does not occur as composition factor of $M$, then
$X=S\oplus M$, and $P$ is indecomposable.
\end{lemma}
\pf Let $\be:P\to S$ be an epimorphism. Then Lemma \ref{homdim}
implies immediately that $P$ is indecomposable with head
$D=\bar{S}/\Jac(\bar{S})$ which is irreducible. Furthermore the
multiplicity of $D$ in $\bar{S}$ is one and $\End_{kG}(\bar{S})\cong k$.

Consider the following diagram, where $\pi:X\to S$ is the canonical
epimorphism: 
$$
\begin{diagram}
&&&&&&\ker\be&&\ker\rho\\
&&&&&&\dTo&\ldTo&\\
&&&&&&P&&\\
&&&&&\ldTo^{\exists\rho}&\dTo_\be&&\\
(0)&\rTo&M&\rTo&X&\rTo^\pi&S&\rTo&(0)
\end{diagram} 
$$
By projectivity of $P$ we find $\rho:P\to X$ such that $\pi\rho =
\be$. Let $U = \ker\be$. Now $\pi\rho(U)=\be(U) =(0)$, hence 
$\rho(U)\subseteq\ker\pi$. Tensoring over $\Or$ by $K$ yields
$$
(1_K\otimes\rho)(KU)\subseteq\ker(1_K\otimes\pi)=KM.
$$ 
We conclude
that either $KU$ and $KM$ have a composition factor in common,
or $KU\leq\ker(1_K\otimes\rho)$, and hence $U\leq\ker\rho$.   
Since $D$ is not a composition factor of $\bar{M}$ we see, using
Brauer reciprocity, that $KM$ and $KP$ do not have a composition factor in
common, therefore the same holds for $KU\leq KP$ and $KM$. 

Thus $U\leq\ker\rho$. We have an induced
map $\hat{\rho}$ from $P/U\cong S$ to $X$ such that
$\rho=\hat{\rho}\alpha$, where $\alpha:P\to P/U$ is the natural
epimorphism. Similarly $\be$ induces an isomorphism $\hat{\be}$ from
$P/\ker{\be}$ to $S$. We have $\pi\hat{\rho}=\hat{\be}$, and therefore
$\pi\hat{\rho}\hat{\be}^{-1}=1_S$ and $\pi$ splits, as desired. 
\epf

Let $G$ be a finite group with split $BN$-pair of 
characteristic $p$. We want to investigate representations of $G$ in
the non-describing characteristic case, hence we assume that the prime
$p$ is  different from $\ell$. The notion of a $BN$-pair
means that we have given a Borel subgroup $B$ and the monomial
subgroup $N$ of $G$. For the definition and details on groups with 
$BN$-pair we refer to \cite{ca1}. Occasionally, $G$ has to be viewed as
a Levi subgroup of a larger finite group with split $BN$-pair $\tG$.
We denote the $BN$-pair of $\tG$ by $\tilde{B}$ and $\tilde{N}$
and say that $\tG$ has a $\tilde{B}\tilde{N}$-pair. In this situation
we often assume that $G$ is a \emph{standard} Levi subgroup of
$\tG$, that is, if $B=UT$ and $\tilde{B} = \tilde{T}\tilde{U}$ is the
Levi decomposition of $B$ and $\tilde{B}$ respectively, then
$T=\tilde{T}$, $B\subseteq\tilde{B}$ and $B\tilde{U}=\tilde{B}$. As a
consequence the monomial subgroup $N$ of $G$, being the normalizer
$\cN_G(T)$, is contained in $\tilde{N}$ and $W=N/T \leq \tilde{N}/T = 
\tilde{W}$, indeed the \emph{Weyl group} $W$ of $G$ is a standard
parabolic subgroup of $\tilde{W}$ generated by a subset of the
standard generators (simple reflections) of $\tilde{W}$.  
 
The set of Levi subgroups of $G$ will be denoted by $\cL_G$
and for $L\in \cL_G$, we denote by $\cL_{G,L}$ the set of Levi
subgroups of $G$ containing $L$. Note that every Levi subgroup of $G$
is conjugate to a standard Levi subgroup.  

For $L\in \cL_G$ and $\Rs$, we define the functor  
\emph{Harish-Chandra induction} $\R_L^G$ as follows: We choose a
parabolic subgroup $P$ of $G$
such that $L$ is a Levi complement of $P$. Let $U$ be the Levi kernel,
that is the unipotent radical, of $P$. Thus $P$ is the semi direct
product $P = U\rtimes L$. If $M$ is an $RL$-lattice we may consider it
as $RP$-lattice with trivial $U$-action. The corresponding functor is
called \emph{inflation} and is denoted by $\Inf_L^P$. The functor
$\R_L^G$ is now defined to be inflation followed by ordinary
induction $\Ind_P^G$. Thus for $M\in\lat_{RL}$ we have:
$$
\R_L^GM = \Ind_P^G\circ\Inf_L^P(M).
$$
Dually we define \emph{Harish-Chandra restriction}$\T_L^G$ to be the 
restriction $\Res^G_P$ from
$RG$-modules to $RP$-modules followed by $\Fix_U$, which is the
functor mapping the $RP$-module $M$ to the $RL$-module $\Fix_U(M) =
\{m\in M\mid mu=m\text{ for all } u\in U\}$. The fact that the order of $U$ 
is invertible
in $R$ implies that $\R_L^G$ and $\T_L^G$ are adjoint functors, and the
collection of functors $\R_L^G$ and their adjoints satisfy
transitivity and a Mackey formula. They were investigated in a more
general context in \cite{difl1} and \cite{didu1}. We point
out that most of our general results in the following sections 
can be generalized
to the situation in \cite{didu1}. Morover, in this paper and
independently in \cite{hole2} it was shown that the functors $\R_L^G$
and their adjoints are independent of the choice of the parabolic
subgroup $P$ of which $L$ is a Levi subgroup, provided the order of
$U$ is invertible in $R$. For fields $R$ of characteristic $0$ this
result is known for a long time and was observed first by Deligne
(comp. e.g. \cite{lusp}). Of course these functors yield corresponding
maps on characters and Brauer characters which again are denoted by
$\R_L^G$ and $\T_L^G$.

Harish-Chandra induction and restriction (or HC-induction and
HC-restriction for short) provide the basic tools for the Harish-Chandra
theory which subdivides the irreducible $RG$-modules for $R = K$ or
$k$ into \emph{Harish-Chandra series} (or HC-series for short). 
First an $RL$-module $M$ is called \emph{cuspidal} if 
$\T_{L'}^L(M)= (0)$ for every proper Levi subgroup $L'$ of $L$. The 
irreducible $RG$-module $V$ belongs to the HC-series $S(G/L,M)_R$
if $V$ is a composition factor of the head (or alternatively of the
socle) of $\R_L^G(M)$. If $R=K$ we usually omit the index $K$. The
subgroup $L$ and the irreducible cuspidal $RL$-module $M$ are unique
up to conjugation in $G$ respectively in the normalizer $\cN_G(L)$ of 
$L$ in $G$. Analogously we define the HC-series $S(G/L,\chi)_R$ for
the irreducible cuspidal (Brauer-) character $\chi$ of $L$. This
partition of the irreducible $RG$-modules and characters into
HC-series is well known for fields $R$ of chracteristic $0$ (see
e.g. \cite{ca1}), and was introduced by Hiss in \cite{hi1} for fields
$R$ of positive characteristic $\ell$. If the $RG$-module $X$ is in
the HC-series $S(G/L,M)$, the Levi subgroup $L$ is called
\emph{semisimple vertex} and the irreducible cuspidal $RL$-module $M$
\emph{semisimple source} of $X$ (see \cite{didu1}). As remarked,
those are unique up to conjugation in $G$ and hence we may (and
usually do) choose as representatives for semisimple vertices the
unique standard Levi subgroups.   

The following two 
lemmas are just formulated for the convenience of the reader. They
can be easily derived from \cite[5.8]{hi1} by Frobenius reciprocity
and Mackey decomposition. Let now $L\leq  M$ be  Levi subgroups 
of $G$. We take an irreducible cuspidal $kL$-module $X$ and an irreducible
$kG$-module $Y$.
\begin{lemma}\label{dualstatement1} 
$Y$ lies in the HC-series $S(G/L,X)_k$ if and only if we find
$Z,Z'\in S(M/L,X)_k$ such that $Y$ is isomorphic to some composition
factor of $\hd(\R_M^GZ)$ and to some composition factor of $\soc({\R_M^GZ'})$.
\end{lemma}

\begin{lemma}\label{dualstatement2} $Y$ lies in the series
  $S(G/L,X)_k$ if and only if there exist $Z,Z'\in S(M/L,X)_k$ such
  that $Z$ is isomorphic to some composition factor of $\hd(\T_M^GY)$ and 
$Z'$ to some composition factor of $\soc(\T_M^GY)$. 
\end{lemma}

\begin{notation}\label{notquot}{\rm 
Let $\cX$ be a set of irreducible characters of $G$ and let $Y$ 
be an $\O G$-lattice affording the character $\chi$. We may write uniquely
$\chi=\psi+\psi'$, where every irreducible constituent of $\psi$ belongs to
$\cX$ and no irreducible constituent of $\psi'$ belongs to
$\cX$. We find pure sublattices $V$ and $V'$ of $Y$ such that $V$ affords
$\psi$ and $V'$ affords $\psi'$. We define $\pi_{\cX}(Y):= Y/V'$, the
\emph{$\cX$-quotient} of $Y$ and $\iota_{\cX}(Y):= V$, the
\emph{$\cX$-sublattice} of $Y$.
If $\cX$ consists of a single character $\phi$ say we write $\pi_{\phi}(Y)$
and $\iota_{\phi}(Y)$ instead of $\pi_{\{\phi\}}(Y)$
and $\iota_{\{\phi\}}(Y)$ respectively. Moreover, we set $\pi_{\cX}(\chi)
=\psi=\iota_{\cX}(\chi)$. We define 
$\pi_{\cX}(\bar{Y}):= \overline{\pi_{\cX}(Y/V')}$ and 
$\iota_{\cX}(\bar{Y}):= \overline{\iota_{\cX}(Y/V')}$. We note that the
latter two depend on $Y$ and not just on $\bar{Y}$.  
Obviously $\pi _{\phi}(Y)$ is zero or affords a multiple of the character
$\phi$. Let $\fC_\cX$ be the full subcategory of the category of finitely 
generated $\O G$-lattices whose objects are those $Y\in\lat_{\O G}$
satisfying $\iota_{\cX}(Y)=Y$.} 
\end{notation}

The construction of $\pi$ is precisely, what is needed to apply the
results from the beginning of this section, as the following results
demonstrates. It follows immediately from the definitions in \ref{notquot}:

\begin{lemma}\label{hompi} Let $Y$ be a projective $\O G$-lattice and
  $\cX$ a set of irreducible characters of $G$. Let $M=\pi_\cX(Y)$ and  
  $\be:Y\to M$ be the natural projection. Then $\be$ satisfies
  Hypothesis \rmref{hypo} and the requirements of Theorem
  \rmref{decmatrix}. In particular the decomposition matrix of the
  $\O$-order $\H=\End_{\O G}(M)$ is part of the $\ell$-modular
decomposition matrix of $G$.
\end{lemma}
  
By construction $\pi $ and $\iota$ are functorial: If we have given  
$Z\in\lat_{\O G}$ with $U=\pi_{\cX}(Z)$ and $U'=\iota_{\cX}(Z)$ and if
$f:Y\to Z$, then $f|_V:V\rightarrow U$ and
$f|_{V'}:V'\rightarrow U'$. Thus we have induced maps $\pi_{\cX}(f):
\pi_{\cX}(Y)\rightarrow \pi_{\cX}(Z)$ and $
\iota_{\cX}(f):\iota_{\cX}(Y):\rightarrow \iota_{\cX}(Z)$. Thus 
$\pi_{\cX}$ and $\iota_{\cX}$ are funtors from $\lat_{\O  G}$ onto 
$\fC_\cX$, indeed both are left inverses of the embedding of $\fC_\cX$
into $\lat_{\O  G}$. Obviously for $\cX=\emptyset$ the functors
$\pi_{\cX}$ and $\iota_{\cX}$ are both the zero functor, taking every
module to the zero module. If $\cX$ contains all irreducible
characters of $G$ we get the identity functors.

\begin{lemma}\label{functprop} Let $Y,Z\in\lat _{\O  G}$. 
Let $\cX$ be a set of irreducible characters of $G$.
If either $Y$ or $Z$ is in $\fC_\cX$ then the maps 
\begin{align*}
\pi_{\cX}: \Hom _{\O G}(Y,Z)&\to \Hom _{\O G}(\pi_{\cX}(Y),\pi_{\cX}(Z))
\intertext{and}
\iota_{\cX}: \Hom _{\O G}(Y,Z)&\to \Hom _{\O G}(\iota_{\cX}(Y),\iota_{\cX}(Z))
\end{align*}
induced by the functors $\pi_{\cX}$ and $\iota_{\cX}$ are
isomorphisms.
\end{lemma}
\pf This is immediate by tensoring by $K$ and comparing dimensions.
\epf

Lemma \ref{functprop} does not carry over to liftable $kG$-modules 
of the form $\bar{Y}$ and $\bar{Z}$ for $Y,Z\in\lat_{\O G}$ 
$kG$-modules, since the canonical embedding $k\otimes _{\O}\Hom_{\O
  G}(Y,Z)\rightarrow\Hom_{k G}(\bar{Y},\bar{Z})$ is not surjective 
in general. However it is a
bijection in case that either $Y$ or $Z$ is a projective $\O G$-lattice (see
e.g. \cite[Theorem 14.7]{la}). 

\begin{lemma}\label{functpropproj} Let $Y,Z\in \lat_{\O G}$ such that 
$\Hom_{kG}(\bar{Y},\bar{Z})$ is $liftable$. Let $\cX$ be a set of
  irreducible characters of $G$ and let $\fC_\cX$ be defined as in
  \rmref{notquot}. Then we have: 
\begin{enumerate} 
\item Let $Z\in \fC_\cX$. Then $\Hom_{k G}(\bar{Y},\bar{Z})\simeq 
\Hom_{kG}(\pi_{\cX}(\bar{Y}),\bar{Z})$.
\item Let $Y\in \fC_\cX$. Then $\Hom_{kG}(\bar{Y},\bar{Z})\simeq \Hom_{k
 G}(\bar{Y}, \iota_{\cX}(\bar{Z}))$.
\end{enumerate}\end{lemma}

\pf By Definition \ref{notquot} we have $\pi_{\cX}(\bar{Y})
=\bar{Y}/\bar{V'}$ for $V'\in\lat_{\O G}$ and $\iota_{\cX}(V')=(0)$. 
Denoting the projection from $\bar{Y}$ onto $\pi_{\cX}(
\bar{Y})$ by $\bar{\pi}$ we have an embedding 
$$
\bar{\pi}^*:\Hom_{kG}(\pi_{\cX}(\bar{Y}),\bar{Z})\rightarrow \Hom_{k
  G}(\bar{Y},\bar{Z}) 
$$
defined by $\bar{\pi}^*(f)=f\circ \bar{\pi}$. By assumption and Lemma
\ref{functprop} we have 
$$
\Hom_{kG}(\bar{Y},\bar{Z}) =k\otimes _{\O}\Hom_{\O G}(Y,Z)\simeq 
k\otimes _{\O}\Hom_{\O G}(\pi_{\cX}(Y),Z)\subseteq \Hom_{kG}(\pi_{\cX}
(\bar{Y}),\bar{Z}),
$$
hence, comparing dimensions, $\bar{\pi}^*$ is an isomorphism.
ii) is proven analogously.
\epf

We shall study now how the functors $\pi$ and $\iota$ behave with respect to
HC-induction and -restriction, where we take in \ref{notquot}
HC-series of characters of $G$ as set $\cX$.

\begin{notation}\label{hcquot} {\rm Let $G$ be a group with split
    $BN$-pair which is a Levi subgroup of a group $ \tG$ with 
split $\tilde{B}\tilde{N}$-pair.
Let $\M\subset\L_{G} $ and let $\C$ be a set of cuspidal
irreducible characters of elements of $\M$. We say the pair $(\M,\C)$ is
\emph{closed under conjugation} in $\tilde{G}$ if $M\in \M$ implies 
$M^x\in \M$
and $\chi\in \C$ implies $\chi^x\in \C$ for all $x\in \tilde{G}$ such that
$M^x\in \L_{G}$. If $M\in \L_{G}$ we set $\M_M:=\M\cap \L_M$, and we
define $\C_M$ to be the set
of those $KM'$-characters $\chi$ with $M'\in \M_M$ and $\chi\in
\C$. Note that
$(\M_M,\C_M)$ is closed under conjugation in $\tilde{G}$ if $(\M,\C)$
is. On the other hand if $(\M,\C)$ is an arbitrary pair for
$G\in\L_{\tG}$, then $(\M,\C)$ is a
pair for $\tG$. However, even if $(\M,\C)$ is closed under
conjugation by $G$, it is not necessarily closed under conjugation by
$\tilde{G}$.

Now let $(\M,\C)$ be fixed. We set for every $M\in \L_G$
$$
\cX_M:=\bigcup_{L\in \M_M,\atop \chi\in Irr(L)\cap
  \C_M}S(M/L,\chi).
$$
Then we define $\pi_{(\M,\C)}^M=\pi_{\cX_M}$. If no ambiguities can arise
we omit the index $(\M,\C)$ and write simply $\pi^M = \pi_{\cX_M}$.}
\end{notation}

\begin{rem}\label{warning} Obviously if $M$ does not contain any
  $L\in\M$ then $\pi_{\cX_M}=\pi_{\emptyset}$ is the zero
  functor. Moreover it should be pointed out that the
  \emph{restriction} $\pi^M$ of $\pi$ to $M\in\L_G$ essentially
  depends upon the choice of $(\M,\C)$, not on the functor $\pi =
  \pi^G$ alone. For instance we choose $x\in G$ such that $L\leq M$ and
  $L^x\leq M$ are not conjugate in $M\in\L_G$, and we take an irreducible
  cuspidal character $\chi$ of $L$. Then $\chi^x$ is an irreducible
  cuspidal character of $L^x$. We set
  $\M=\{L\},\;\C=\{\chi\},\;\M^x=\{L^x\}$ and $\C^x=\{\chi^x\}$. Then
  obviously
$$
\pi^G_{(\M,\C)}=\pi^G_{(\M^x,\C^x)}
$$ 
but
$$
\pi^M_{(\M,\C)}\not\cong\pi^M_{(\M^x,\C^x)}.
$$
\end{rem}

To apply Results \ref{notquot} and \ref{functpropproj} we need 
the following lemma which
follows immediately from Mackey decomposition (compare \cite[1.14]
{difl1}) and from transitivity of HC-induction:

\begin{lemma}\label{hcserprel} Let $\psi$ be a $KG$-character.
\begin{enumerate}
\item Let $M\in \L_{G,L}$ and assume that all constituents of
  $\psi$ belong to the HC-series $S(G/L,\chi)$, where $\chi$ is an 
irreducible cuspidal $KL$-character. Then every constituent of 
$\T^G_M\psi$ belongs to a series $S(M/L^x,\chi^x)$ for some $x\in G$ 
such that $L^x\leq M$.
\item Suppose $G\in \L_{\tG}$ and assume that every irreducible
constituent of $\psi$ belongs to some  HC-series $S(G/L^x,\chi^x)$, where
$x\in \tilde{G}$ such that $L^x\leq G$. Then $\R_G^{\tG}\psi $
belongs to $S(\tG/L,\chi)$.
\end{enumerate}\end{lemma}
\pf To prove (i) we may assume that $\psi$ is irreducible. Then $\psi$ is
constituent of $\R_L^G\chi$ by assumption and hence $\T^G_M\psi$ is a 
summand of
$$
\T^G_M\R_L^G\chi=\sum_x\R_{L^x\cap M}^M\T^{L^x}_{L^x\cap M}\chi^x,
$$
where $x$ runs through a set of double coset representatives of parabolic
subgroups having $M$ and $L$ as Levi complements, by
\cite[1.14]{difl1}. But
these double coset representatives can be chosen in $N$, 
(see \cite[2.1]{ca1}).
Moreover $\T^{L^x}_{L^x\cap M}\chi^x=(0)$ unless $L^x\cap M=L^x$, that is
$L^x\leq M$ by cuspidality of $\chi$. Thus the result follows.
Part (ii) is shown similarly, observing that for $x\in N$ we have
$\R_{L^x}^{\tG}\simeq \R_L^{\tG}$ by \cite{hole2} respectively 
\cite[5.2]{didu}.
\epf

Note that in the lemma above we may take $M$ and $L$ to be
standard Levi subgroups of $G$, and $G$ to be a standard Levi subgroup
of $\tG$ (by adjusting $G$, then $M$ and then $L$ by conjugating by an
element of $\tG\;,G$ and $M$ respectively). In this case we may take
$x$ to be an element of $N$ or in (ii) of $\tilde{N}$. Now 
Lemma \ref{hcserprel} tells us that under HC-restriction we might get 
terms with semi simple HC-vertex in $M^x$, where $x\in N$ such that 
$L^x\leq M$.
We therefore define 

\begin{ddef}\label{nlm} {\rm Let $M\in \L_{G,L}$. Then we set
$$
\NLM =\{x\in G|L^x\leq M\}.
$$
The letter $N$ stands here for the monomial group of $G$, since we can
always assume that $L$ is standard in $M$ and $M$ in $G$, and then
choose $x$ to be an element of $N$. Thus for instance 
$$
\tilNLG=\{x\in\tG\mid L^x\subseteq G\}.
$$}
\end{ddef}

The following Corollary follows now easily from Lemma \ref{hcserprel}:

\begin{cor}\label{commute} Suppose that $(\M,\C)$ is closed under conjugation
  in $G$ and let $M\in \L_G$. Then we have isomorphisms of functors
\begin{enumerate}
\item For $(\M,\C)$ as in \rmref{hcquot} and $M\in\L_G$ we have 
\begin{align*}
\R_M^G\pi^M&\simeq \pi^G\R_M^G\\
\R_M^G\iota^M&\simeq \iota^G\R_M^G.
\end{align*}
\item 
\begin{align*}
\T^G_M\pi^G&\simeq \pi^M\T^G_M\\
\T^G_M\iota^G&\simeq \iota^M\T^G_M.
\end{align*}
\end{enumerate}\end{cor}

\pf First note that it suffices to check that the functors agree on
characters, since $\R_M^G$ and $\T_M^G$ preserve purity of
sublattices. We show the second part of Corollary, the proof of the
first being similar. So let $\psi$ be an irreducible 
$KG$-character in $S(G/L,\chi)$, and $L\in \M$. If $\psi\in \C$ then 
$\pi^G(\psi)=\psi$
and $\pi^G(\psi)=(0)$ otherwise. If $L^x\not\leq M$ for all $x\in
G$, then $\T^G_M\psi=(0)$ by Mackey decomposition and then 
$$
0=\T^G_M\pi^G(\psi)=\pi^M(\T^G_M\psi)
$$
in all cases. Suppose $L^x\leq M$. Since
$S(G/L,\chi)=S(G/L^x,\chi^x)$ we
may assume $x=1$, that is $L\leq M$. Let $\chi\in \C$ then 
$\chi^x\in \C_M$ for all $x\in\NLM$ and hence 
$$
\pi^M(\T^G_M\psi)=\T^G_M\psi=\T^G_M\pi^G(\psi).
$$
If $\chi\notin\C$ then $\chi^x\not\in \C$ for all $x\in \NLM$ 
and we get
$$
0=\pi^M(\T^G_M\psi)=\T^G_M(0)=\T^G_M\pi^G(\psi).
$$
The other isomorphisms of functors are shown similarly 
\epf

\begin{cor}\label{projhom} Suppose that $(\M,\C)$ is closed under 
conjugation in $G$ and let $M\in \L_G$. Suppose $Y\in\lat_M$ is
projective. Set $X=\pi^M(Y)$ and consider 
$$
\be:\R_M^GY\to R_M^GX,
$$
where $\be = \R_M^G\psi$ setting $\psi:Y\to X$ to be the natural
projection. Then $\be$ satisfies the requirements of Theorem
\rmref{decmatrix}. In particular the decomposition matrix of
$\H=\End_{\O G}(\R_M^G(X))$ is part of the $\ell$-modular decomposition
matrix of $G$.
\end{cor}

\pf By Corollary \ref{commute}
$$
\pi ^{G}(\R_M^{G}Y)=\R_M^{G}\pi^M(Y)=\R_M^{G}X.
$$
Thus the result follows from \ref{hompi} applied to the projective $\O
G$-module $\R_M^GY$. \epf

\section{Projective restriction systems}\mlabel{projrestsys}

We introduce now \emph{projective restriction systems} for finite
groups with split $BN$-pairs. First we describe the set up which we
shall be dealing with:

\begin{notation}\mlabel{notrestsys} {\rm Let $G$ be a finite group with
split $BN$-pair. We fix $L \in \L_G$ and an irreducible cuspidal 
$KL$-character $\chi_L$. We assume that $\chi_L$ occurs 
with multiplicity one as an irreducible constituent of the character
afforded by some projective $\O L$-lattice $Y_L$. Thus we may apply
Lemmas \ref{homdim} and \ref{inprel2}.

If $G\in \L _{\tG}$ for some group $\tG$ with
$\tilde{B}\tilde{N}$-pair such that $B=\tilde{B}\cap G$ and
$N=\tilde{N}\cap G$, we take in \ref{hcquot}  
$$
\M :=\{L^x\mid x\in \tilNLG\}
\quad\text{and}\quad\C:=\{\chi_L^x\mid x\in \tilNLG\}.
$$
For $M\in \L_G$ we write 
$$
\pi^M:=\pi ^M_{(\M,\C)}.
$$
In particular for $M=G$ we omit the superscript $G$ and write
$\pi=\pi^G$. A similar notation is used for $\iota$. 

Finally we set 
$$
X_L:=\pi^L(Y_L),\quad X_L':=\iota^L(Y_L),\quad\text{and}\quad 
D_L:=\bar{X}_L/\Jac(\bar{X}_L).
$$} 
\end{notation}

\begin{rem}\mlabel{xldifxl} 
By Lemma \rmref{inprel2}, $D_L\cong\hd(\bar{X}_L)$ is irreducible and is
  the head of the projective cover $\bar{Y}_1$ of $\bar{X}_L$ which is
  isomorphic to a direct summand of $\bar{Y}_L$. If $Y_1$ denotes
  its lift to an $\O L$-sublattice of $Y_L$, which exists by
  projectivity of $Y_L$, then $\iota^L(Y_L)\subseteq Y_1$.
Since $kL$ is a symmetric algebra, 
the head $D_L$ of $\bar{Y}_1$ and its socle are isomorphic. Thus
the lattice $X_L'$ is different from $X_L$ unless
$\bar{X}_L$ is irreducible, as $D_L$ is the head of  $\bar{X}_L$   and
the socle of $\bar{X}_L'$. Obviously we have
$$
\Hom_{kL}(\bar{X}_L,\bar{X}_L')\simeq k.
$$
\end{rem}

\begin{ddef}\mlabel{defrestsys} {\rm Keep the notation introduced in
\ref{notrestsys}. In view of Lemma \ref{homdim} we may assume that
$Y_L$ is indecomposable and hence is the projective cover of
$X_L$. Note that this implies that $\bar{Y}_L$ is an principal
indecomposable $kL$-module corresponding to the irreducible
module $D_L$. In particular $D_L$ is the head
$\bar{Y}_L/\Jac(\bar{Y}_L)$ of $\bar{Y}_L$. Suppose in addition that 
$\chi_G$ is an irreducible constituent of multiplicity one in
$\R_L^G\chi_L$ such that the following holds:
\begin{enumerate}
\item There exists a projective $\O G$-lattice $Y_G$ such that
$\chi_G$ is the character of 
$$
\pi^G(Y_G)=:X_G.
$$
Again we may (and do) assume that $Y_G$ is indecomposable.
\item For $M\in \L_{G,L}$ let $\T^G_M\chi_G=\psi_1+\ldots +\psi_n$
be a decomposition of $\T^G_M\chi_G$ into irreducible constituents. Let
$X_i$ be an $\O M$-lattice affording $\psi _i$ for $i=1,\ldots ,
n$. Then for $i\not= j$, the $kM$-modules $\bar{X}_i$ and $\bar{X}_j$
have no composition factor in common.
\end{enumerate}
We set $X_M:=\pi_{(L,\chi_L)}(\T^G_MX_G)$ for every $M\in \L_{G,L}$
and $X_M':=\iota_{(M,\chi_M)}(\R_L^MX_L)$. We
denote the characters of $X_M$ by $\chi_M$ and the projective
cover of $X_M$ by $Y_M$.

The \emph{projective restriction system} to $(L,\chi_L)$ then consists
of the data $\{X_M,Y_M,X_M'\mid M\in \L_{G,L}\}$ and is denoted by
$\PR(X_G,Y_L)$.}
\end{ddef}

\begin{rem}\mlabel{consistent}
Our notation is 
consistent, since part {\rm (i)} of \rmref{defrestsys} implies 
$$
\pi (Y_G)=\pi^G_{(\M,\C)}(Y_G)=\pi_{(L,\chi_L)}(Y_G)=X_G.
$$ 
Hence $X_G=\pi_{(L,\chi_L)}(\T^G_GX_G) $ since $X_G$ is irreducible and
$\chi_G\in S(G/L,\chi_L)$ by assumption. Similarly $Y_L$ is the
projective cover of $\pi_{(L,\chi_L)}\T^G_LX_G = X_L$, since 
$$
\Hom_{\O L}(\T_L^GX_G,X_L)\cong\Hom_{\O G}(X_G,\R_L^GX_L)\cong\O
$$
by Frobenius reciprocity. Moreover, if $x\in \cN_G(L)$ is such that 
the conjugate character $\chi_L^x$ is different from $\chi_L$, then
$\chi_L^x$ cannot be a constituent of $Y_L$. With other words we have
$\pi^L(Y_L)=X_L$. Otherwise the simple head
of $\bar{Y}_L$ is composition factor of $\bar{X}_L$ as well as of
$\bar{X}_L^x$, since $Y_L$ is indecomposable. But
$X^x_L$ affords $\chi_L^x$, and it is constituent of $\T^G_{L}X_G$,
hence $\bar{X}_L^x$ has by {\rm (ii)} of \rmref{defrestsys} no 
composition factor with $\bar{X_L}$ in common. 
We also note that part {\rm (ii)} is particularly satisfied if the
summands of $\T^G_M\chi_G$ are in pairwise different
$\ell$-modular blocks of $\O M$. This will be the fact which we prove in
concrete applications.
\end{rem}

For the remainder of this section we fix a projective restriction system
$\PR(X_G,Y_L)$. Throughout $M\in\L_{G,L}$. The next lemma is trivial: 

\begin{lemma}\mlabel{conjprojrestsys} Let $x\in\tG$, then 
the conjugate projective restriction system 
$$
\PR^x(X_G,Y_L)=\PR(X_G^x,Y_L^x)
$$
is the projective restriction system given as 
$\{X_M^x,Y_M^x,X_M^{x \prime}\mid M\in\L_{G,L}\}$.
\end{lemma}

Observe that $X_M^x$ is the conjugate module $M^x$-module $(X_M)^x$.
We have to distinguish this from certain other modules coming up in the
conjugate restriction system:

\begin{notation}\mlabel{moreprojnot} Let $x\in\tG$ and
  $M\in\L_{G^x,L^x}$. Then the $M$-module
$$
X_{M,x} = \pi_{(L^x,\chi^x_L)}(\T^G_MX_G^x)
$$ 
and its projective cover $Y_{M,x}$ belong to the conjugate restriction
system $\PR^x(X_G,Y_L)$. In particular, if  $M\in\L_{G,L}$ and $x\in G$
is such that $L^x\subseteq M$, then 
$$
\PR^x(X_G,Y_L) = \PR(X_G,Y_L^x)
$$
and 
$$
X_{M,x} = (X_{M^{x^{-1}}})^x = \pi_{(L^x,\chi^x_L)}(\T^G_MX_G).
$$
We denote the related terms similarly: For example $\chi_{M,x}$
denotes the character afforded by $X_{M,x}$.
\end{notation}

In the following
we shall use frequently the Mackey formula for HC-induction and
restriction in connection with the property of $X_L$ to be
cuspidal. In the Mackey formula we deal with double coset
representatives $\cD_{LM}$ of parabolic subgroups containing $L$
respectively $M\in\L_{G,L}$ as Levi complement. If $L$ and $M$ are
standard Levi subgroups we may choose $\cD_{LM}$ to be
\emph{distinguished} by general theory, that is it is contained in
the monomial subgroup $N$ of $G$ as set of preimages of the
distinguished double coset representatives of the corresponding
standard parabolic subgroups of the Weyl group $W=N/T$ of $G$, (for
details see e.g. \cite{ca1}). This choice then ensures, that the 
intersections $M\cap L^x$ with $x\in\cD_{LM}$ are again standard Levi 
subgroups of $G$. Since $X_L$ is cuspidal by assumption, all terms in 
the Mackey formula
\begin{equation}\mlabel{mackey}
\T^G_M\R^G_LX_L = \bigoplus_{x\in\cD_{LM}}\R_{M\cap L^x}^M\T_{M\cap
  L^x}^{L^x}X_L^x
\end{equation}
with $L^x\nleq M$ are zero, that is we have to sum only over
$x\in\NLM\cap\cD_{LM}$. But $x\in G$ is contained in $\NLM$
precisely if $LxM = xM$. If in addition $x\in\cD_{LM}$, then $L^x =
L^x\cap M\in\L_G$. Obviously $M$ acts on
$\NLM$ by right translation, where the orbit of $x\in\NLM$ under this
action is the double coset $xM =LxM$. 

Observe that $\cN_G(L)\subseteq \NLM$, and that this
group acts on $\NLM$ by left translation. We take $\cNNML$ to be a set
of representatives of $\cN_G(X_L)\times M$-orbits on $\NLM$
(containing $1\in\NLM$),
where $\cN_G(X_L)$ is the stabilizer of $X_L$ in $\cN_G(L)$. Note that
each such orbit contains elements of $\cD_{LM}$ and we choose therefore
$\cNNML$ to be a subset of our selected set of double coset
representatives $\cD_{LM}$. If we need to
emphasize the dependence of the chosen object $X_L$ we write
$\cNNML(X_L)$, and we define similar representatives for other
$L$-modules as for example $\cNNML(D_L)$ etc. 

It follows now that $\{L^x\mid x\in \cNNML\}$ is a 
set of Levi subgroups of $M$ which for $x\neq 1$ are  
conjugate to $L$ in $G$ but not in $M$. Moreover $L^x=L^y$ with 
$x,y\in \cNNML$ and $x\not=y$, implies $X_L^x\not\simeq X_L^y$.

\begin{thm}\mlabel{structhcrest} Let $M\in \L_{G,L}$ and
  $x\in\cNNML$. Then the $\O M$-lattice $X_{M,x}$ affords the
  character $\chi_{M,x}$ which is irreducible and occurs with multiplicity
  one as constituent of the character afforded by $\T_M^G(Y_G)$. Moreover 
$$
\T^G_MX_G=\bigoplus\limits_{x\in\cNNML}X_{M,x}.
$$
\end{thm}

\pf The character $\chi_G$ of $X_G$ is contained in
the HC-series $S(G/L,\chi_L)$. Hence by Lemma \ref{hcserprel} every 
irreducible constitutent of $\T^G_M\chi_G$ belongs to some 
$S(M/L^x,\chi^x_L)$ for some $x\in\NLM$. We may assume that
$x\in\cNNML$. By assumption $\chi_G$ occurs with multiplicity one in 
$\R^G_L\chi_L$, hence by Lemma \ref{conjprojrestsys} $\chi_G$ occurs
with multiplicity one in $\R_{L^x}^G\chi_L^x$ too. Therefore by
Frobenius reciprocity and transitivity of HC-restriction
$$
K\cong \Hom_{KL^x}(KX_L^x,\T_{L^x}^G KX_G)\cong
\Hom_{KM}(\R_{L^x}^MKX_L^x,\T^G_MKX_G),
$$
and consequently $\T_M^G\chi_G$ has precisely one constituent in
$S(M/L^x,\chi^x_L)$ and this is by definition given as 
$$
\pi_{L^x,\chi_L^x}(\T_M^G\chi_G) = \chi_{M,x}.
$$
In particular $\chi_{M,x}$ is an irreducible character and, being in 
different HC-series of $M$, we have $\chi_{M,x}\not=\chi_{M,y}$ for 
$x\not=y$ in $\cNNML$. We have shown so far that
$$
\T^G_M\chi_G=\sum\limits_{x\in\cNNML}\chi_{M,x}.
$$
Being an $\O M$-lattice in an irreducible $KM$-module, $X_{M,x}$ is
indecomposable. Note that $X_G$ is an epimorphic image of $Y_G$ and
hence $X_{M,x}$, being an epimorphic image of $\T^G_MX_G$, is an epimorphic
image of the projective $\O M$-lattice $\T^G_MY_G$. Hence the
projective cover $Y_{M,x}$ of $X_{M,x}$ is a direct summand of 
$\T^G_MY_G$. In particular $\chi_{M,x}$ occurs as constituent of the
character of $\T_M^GY_G$. Remark \ref{consistent} in conjunction with Lemma
\ref{conjprojrestsys} implies
$$
X_G=\pi^G(Y_G)=\pi_{(L,\chi_L)}(Y_G) = \pi_{(L^x,\chi_L^x)}(Y_G).
$$ 
By Frobenius reciprocity, Corollary \ref{commute} and Lemma
\ref{functprop} we have
\begin{equation*}
\begin{split}
\Hom_{\O M}(\T^G_MY_G,X_{M,x})
& \cong\Hom_{\O G}(Y_G,\R_M^GX_{M,x})\\
& \cong\Hom_{\O G}(Y_G,\R_M^G\pi^MX_{M,x})\\
& \cong\Hom_{\O G}(Y_G,\pi^G\R_M^GX_{M,x})\\
& \cong\Hom_{\O G}(Y_G,\R_M^GX_{M,x})\\
& \cong\Hom_{\O G}(Y_G,\pi^G\R_M^G(X_{M^{x^{-1}}})^x)\\
& \cong\Hom_{\O G}(Y_G,\pi^G\R_{M^{x^{-1}}}^GX_{M^{x^{-1}}})\\
& \cong\Hom_{\O G}(Y_G,\pi^G\R_M^GX_M)\\
& \cong\Hom_{\O G}(\pi^GY_G,\pi^G\R_M^GX_M)\\
& \cong\Hom_{\O G}(X_G,\pi^G\R_M^GX_M)\\
& \leq 1
\end{split}\end{equation*}
We conclude that the multiplicity of
$\chi_{M,x}$ in the character of $Y_{M,x}$ is precisely one. But now part
(ii) of Definiton \ref{defrestsys} states exactly that the
hypothesis of Lemma \ref{inprel2} are satisfied for $\T^G_MX_G$, $X_{M,x}$
and the projective cover $Y_{M,x}$ of $X_{M,x}$. Thus $X_{M,x}$ is a
direct summand of $\T^G_MX_G$ for all $x\in\cNNML$ and the theorem follows. 
\epf

Note that Theorem \ref{structhcrest} states in particular that $KX_M$
is irreducible. Result \ref{homdim} now implies:
 
\begin{cor}\mlabel{multdm} Let $M\in \L_{G,L}$ and $x\in\cNNML$. Then
  $D_{M,x}:=\bar{X}_{M,x}/\Jac(\bar{X}_{M,x})$ is an irreducible $kM$-module
  occurring with multiplicity one as composition factor of
  $\bar{X}_{M,x}$. Moreover 
$$
\dim_k(\End_{kM}(\bar{X}_M))=1,
$$
and $\bar{Y}_{M,x}$ is the projective cover of $D_{M,x}$. In
particular $\bar{Y}_{M,x}$ and hence $Y_{M,x}$ is indecomposable.
\end{cor}

Note that 
\begin{equation}\label{2.9a}
D_{M,x} = \hd(\bar{X}_{M,x})=\hd((\bar{X}_{M^{x^{-1}}})^x)=
(D_{M^{x^{-1}}})^x.
\end{equation}

As a consequence of Theorem \ref{structhcrest} and part (ii) of
Definition \ref{defrestsys} we get the following two results:

\begin{cor}\mlabel{hcrestx} Let $M,M'\in \L_{G,L}$ with $M'\leq M$.
\begin{enumerate}
\item Then $\pi^{M'}_{(L,\chi_L)}(\T^M_{M'}X_M)=X_{M'}$. In fact 
$$
\T^M_{M'}X_M=\bigoplus\limits_{y\in\cnorm{M'}{\,L}{M}}X_{M',y}.
$$
\item Let $V$ be an indecomposable projective $kM'$-module such that
$V/\Jac(V)$ is a composition factor of $\bar{X}_{M'}$. Then
\begin{align*}
\Hom_{kM'}(\T^M_{M'}\bar{X}_M,V)&\simeq \Hom _{kM'}(\bar{X}_{M'},V)
\intertext{and}
\Hom_{k M'}(V,\T^M_{M'}\bar{X}_M)&\simeq \Hom _{kM'}(V,\bar{X}_{M'}).
\end{align*}
\end{enumerate}
\end{cor}

\begin{cor}\mlabel{projrestsyssubgrp} Let $M\in \L_{G,L}$. Then 
$$
\PR_M:=\PR(X_M,Y_L)=\{X_{M'},Y_{M'},X_{M'}'\mid M'\in \L_{M,L}\}
$$
is a projective restriction system to $(L,\chi_L)$.
\end{cor}

We now investigate the projective cover $Y_G$ of $X_G$:

\begin{thm}\mlabel{projdec} Let $M\in \L_{G,L}$. Then 
$$
\T^G_MY_G=\bigoplus \limits_{x\in\cNNML}Y_{M,x}\oplus
Z,
$$
where $Z$ is a projective $\O M$-lattice with
$\pi^M(Z)=(0)$. Moreover, 
$$
\pi_{(L^x,\chi^x_L)}(\T^G_MY_G)=X_{M,x},
$$
and
\begin{equation*}
\pi_{(L^x,\chi^x_L)}(Y_{M,y}) =
\begin{cases}
\pi_{(L^x,\chi^x_L)}(\T^G_MY_G)
=X_{M,x} &\quad\text{if $x=y$}\\ 
(0)&\quad\text{otherwise},
\end{cases}
\end{equation*}
where $x,y\in\cNNML$. Thus in particular
$$
\pi^M_{(\M,\C)}(Y_{M,x}) = \pi_{(L^x,\chi^x_L)}(Y_{M,x}) = X_{M,x}.
$$
\end{thm}
\pf Since $\T^G_MX_G$ is an epimorphic image of $\T^G_MY_G$,
the latter contains the projective cover of the former as a direct
summand, hence by Theorem \ref{structhcrest}
$$
\T^G_MY_G=\bigoplus_{x\in\cNNML}Y_{M,x}\oplus Z,
$$
for some projective $\O M$-lattice $Z$. Now by Results \ref{commute} 
and \ref{structhcrest}
\begin{equation*}
\begin{split}
\bigoplus_{x\in\cNNML}X_{M,x}\oplus \pi^M(Z)&
= \bigoplus_{x\in\cNNML}\pi^M(Y_{M,x})\oplus\pi^M(Z)\\
&=\pi^M(\T^G_MY_G)=\T^G_MX_G\\
&=\bigoplus \limits_{x\in\cNNML}X_{M,x}=
\bigoplus \limits_{x\in\cNNML}\pi_{(L^x,\chi^x_L)}(\T^G_MX_G).
\end{split}\end{equation*}
Hence $\pi^M(Z)=(0)$ and 
$$
X_{M,x}=\pi_{(L^x,\chi^x_L)}(\T^G_MY_G)=
\pi_{(L^x,\chi^x_L)}(Y_{M,x}).
$$
\epf

\begin{cor}\mlabel{commut} Let $M\in \L_{G,L}$. Then 
$$
\pi_{(L^x,\chi_L^x)}(\R_M^{\G}Y_M)=\R_M^{\G}X_M,
$$
for every $x\in\tG$.
\end{cor}
\pf  We only need to observe that on $\O\G$-lattices 
$\pi^{\G}=\pi^{\G}_{(\M,\C)}=\pi_{(L^x,\chi_L^x)}$ for every $x\in
\tG$, since $L^x$ and $L$ are conjugate in $\G$. Thus by Lemma
\ref{commute}
\begin{align*}
\pi_{(L,\chi_L)}(\R_M^{\G}Y_M)&=\pi^{\G}(\R_M^{\G}Y_M)\\
&=\R_M^{\G}\pi^MY_M\\
&=\R_M^{\G}X_M
\end{align*}
\epf

\begin{notation}\mlabel{irredmod} {\rm Let $L\in \L_{G,L}$. In 
Theorem \ref{structhcrest} we
have seen that $\bar{X}_M$ has a unique maximal submodule, hence
$D_M:= \bar{X}_M/\Jac(\bar{X}_M)$ is irreducible. Now $D_M$ can be
cuspidal or not. We denote the subset of $\L_{G,L}$ of all $M$ such
that $D_M$ is cuspidal by $\L^c(X_G,Y_L)$. Note that $L\in
\L^c(X_G,Y_L)$ since $\chi_L$ and hence $D_L=\bar{X}_L/\Jac(\bar{X}_L)$ is
cuspidal. Note too, that for $x\in G$ the set $\L^c(X_G,Y_L^x)$
consists precisely of
the subgroups $M^x\in\L_{G,L^x}$ with $M\in\L^c(X_G,Y_L)$. In
particular, $M\in\L^c(X_G,Y_L)$ does not imply
$M\in\L^c(X_G,Y_L^x)$ in general, hence if $D_M$ is cuspidal and
$x\in\NLM$ the $kM$-module $D_{M,x}$ is not necessarily
cuspidal. However the irreducible $kM^x$-module $D_M^x=
(D_M)^x$ is cuspidal in this case.} 
\end{notation}

We fix now $M\in\L^c(X_G,Y_L)$ and $G'\in\L_{G,M}$. Note that
$G'\in\L_{G,L}$, and that $\NMGd\subseteq\NLGd$. Take $x\in\NMGd$. The
$\cN_G(X_L)\times G'$-orbit $\cN_G(X_L)xG'$ of $x$ in 
$\NLGd$ is denoted by $[x]_{G',X_L}$ and we may assume that 
$x\in\cNNGdL$. Thus $\bar{X}_{G',x}$ is one of the indecomposable
direct summands of $\T_{G'}^G\bar{X}_G$ by Theorem \ref{structhcrest}. 

\begin{thm}\mlabel{factor} Keep the notation introduced above. Then
  there is precisely one composition factor $D_{G',M^x}$ of 
$\T_{G'}^G\bar{X}_G$ in the HC-series $S(G'/M^x,D_M^x)_k$ and it
is in fact a composition factor of the direct summand $\bar{X}_{G',x}$ of 
$\T_{G'}^G\bar{X}_G$. Moreover, if $V\leq \T_{G'}^G\bar{X}_{G,x}$,
then $D_{G',M^x}$ is a composition factor of $V$ if and only if
$D_M^x$ is a composition factor of $T^{G'}_{M^x}V$. 
\end{thm}

\pf  Using Frobenius reciprocity, part (ii) of Definition
\ref{defrestsys}, and Theorem \ref{structhcrest}, we have 
\begin{align*}
\Hom _{kG'}(\T^G_{G'}\bar{X}_G, \R_{M^x}^{G'}D_M^x)
&=\Hom_{kM^x}(\T^G_{M^x}\bar{X}_G^x,D_M^x)\\
&=\bigoplus_{y\in\cnorm{M^x}{L^x}{G}}
\Hom _{kM^x}(\bar{X}_{M^x,y}^x,D_M^x)\\
&=\Hom_{kM^x}(\bar{X}_M^x,D_M^x)\simeq k,
\end{align*}
since in the last direct sum above the only nonzero term occurs for
$y=1$, again by part (ii) of Definition \ref{defrestsys}, (applied
to the projective restriction system $\PR^x(X_G,Y_L)$).
We conclude that $\T^G_{G'}\bar{X}_G$ and the socle
of $\R_{M^x}^{G'}D_M^x$ have a common composition factor. But this 
is a composition factor $D_{G',M^x}$ of $\T^G_{G'}\bar{X}_G$ in
$S(G'/M^x,D_M^x)_k$.

Now let $V\leq\T^G_{G'}\bar{X}_G^x$. Then the embedding
$V\to\T^G_{G'}\bar{X}_G^x$ induces an embedding of 
$\Hom_{kG'}(\R^{G'}_{M^x}\bar{Y}_M^x,V)$ into
$\Hom_{kG'}(\R^{G'}_{M^x}\bar{Y}_M^x,\T^G_{G'}\bar{X}_G^x)$, and
we have by Theorem \ref{structhcrest}, using \ref{defrestsys},
Frobenius reciprocity, and in addition the fact that $Y_M$ is projective:

\begin{equation}\mlabel{factorequ}\begin{split}
\Hom_{kM^x}(\bar{Y}_M^x,\T^{G'}_{M^x}V)
&\simeq\Hom_{kG'}(\R^{G'}_{M^x}\bar{Y}_M^x,V)\\
&\leq\Hom_{kG'}(\R^{G'}_{M^x}\bar{Y}_M^x,\T^G_{G'}\bar{X}_G^x)\\
&\simeq \Hom _{kM^x}(\bar{Y}_M^x,\T^G_{M^x}\bar{X}_G^x)\\
&\simeq \bigoplus_{y\in\cnorm{M^x}{L^x}{G}}\Hom _{kM^x}
(\bar{Y}_{M^x}^x,\bar{X}^x_{M^x,y})\\
&\simeq \Hom _{kM^x}(\bar{Y}_M^x,\bar{X}_M^x)\\
&\simeq k\otimes _{\O}\Hom_{\O M^x}(Y_M^x,X_M^x)
\simeq k\otimes _{\O}\O\simeq k.
\end{split}\end{equation}

Here we wrote $\bar{X}^x_{M^x,y}$ for the construction
\ref{moreprojnot} applied to the conjugate module $\bar{X}_M^x$.  

The first isomorphism in Equation \ref{factorequ} implies that
there is a non-trivial homomorphism from
$\R^{G'}_{M^x}\bar{Y}_M^x$ to
$V$ if and only if $D_M^x$ is a composition factor of
$\T^{G'}_{M^x}V$. But in this case
$$
\Hom_{kG'}(\R^{G'}_{M^x}\bar{Y}_M^x,V)\simeq
\Hom_{kG'}(\R^{G'}_{M^x}\bar{Y}_M^x,\T^G_{G'}\bar{X}_G).
$$
We observe that $\bar{Y}_M^x$ is the projective cover of
$D_M^x$, hence the projective cover  of 
$\R_{M^x}^{G'}D_M^x$ is a direct summand of 
$\R^{G'}_{M^x}\bar{Y}_M^x$ which in turn contains the projective
cover $\bar{Y}_{G',M^x}$ of $D_{G',M^x}$, since this is 
contained in the HC-series $S(G'/M^x,D_M^x)_k$ and hence is 
an irreducible summand of $\hd(\R_{M^x}^{G'}D_M^x)$. Thus a
nonzero homomorphism from 
$\R^{G'}_{M^x}\bar{Y}_M^x$ into $\T^G_{G'}\bar{X}_G^x$ (which
exists and is unique up to scalar multiplies by \ref{factorequ}), is
nonzero on $\bar{Y}_{G',M^x}$. From this the last assertion of the theorem
follows easily.

It remains to show that $D_{G',M^x}$ is a composition factor of
$\bar{X}_{G',x}$. To do this we apply Theorem \ref{structhcrest} and
the fact that 
homomorphisms from projective modules over $k$ are liftable to $\O$
and obtain:
$$
\O\cong\Hom_{\O G'}(\R^{G'}_{M^x}Y_M^x,T^G_{G'}X_G)\cong
\bigoplus_{y\in\cnorm{G'}{L}{G}}
\Hom_{\O G'}(\R^{G'}_{M^x}Y_M^x,X_{G',y}),
$$
hence there is precisely one $y\in\cnorm{G'}{L}{G}$ such that
\begin{equation}\mlabel{factoreq2}
\Hom_{\O G'}(\R^{G'}_{M^x}Y_M^x,X_{G',y})\neq (0).
\end{equation}
We tensor by $K$ and use \ref{hcquot}, Corollary \ref{commute}, Lemma
\ref{functprop} and Frobenius reciprocity to get:
\begin{equation}\mlabel{factoreq3}\begin{split}
\Hom_{KG'}(\R^{G'}_{M^x}KY_M^x,KX_{G',y})&
=\Hom_{KG'}(\pi^{G'}\R^{G'}_{M^x}KY_M^x,KX_{G',y})\\
&=\Hom_{KG'}(\R^{G'}_{M^x}\pi^{M^x}KY_M^x,KX_{G',y})\\
&=\Hom_{KG'}(\R^{G'}_{M^x}KX_M^x,KX_{G',y})\\
&=\Hom_{KM^x}(KX_M^x,\T^{G'}_{M^x}KX_{G',y})\\
&\leq\Hom_{KM^x}(\R_{L^x}^{M^x}KX_L^x,\T^{G'}_{M^x}KX_{G',y})\\
&=\Hom_{KL^x}(KX_L^x,T^{G'}_{L^x}KX_{G',y}).
\end{split}\end{equation} 
Now $\pi_{(L^x,\chi_L^x)}(X_{G',y})\neq (0)$ if and only if $y=x$,
since $X_{G',y}=\pi_{(L^x,\chi_L^x)}(\T^G_{G'}X_G)$. Thus the first 
$\Hom$-set in Formula \ref{factoreq3} is zero for $y\neq x$, and hence
nonzero for $x=y$ by Formula \ref{factoreq2}. We conclude that 
$$
\Hom_{kM^x}(\bar{Y}_M^x,\T^{G'}_{M^x}\bar{X}_{G',x})=
k\otimes_\O\Hom_{\O G'}(\R^{G'}_{M^x}Y_M^x,X_{G',x})\neq (0). 
$$
and therefore $D_{G',M^x}$ is a composition factor of
$\bar{X}_{G',x}$ by the first part of the proof applied to
$V=\bar{X}_{G',x}$.
\epf

The previous theorem yields a unique composition factor $D_{G',M^x}$ of
 $\bar{X}_{G',x}$ from the HC-series $S(G'/M^x,D_{M^x}^x)$, for every 
$x\in\cNNGdL$ such that $[x]_{G',L}\cap\NMGd$ is not empty. 

Let $x,y\in\NMGd$. Then $S(G'/M^x,D_{M^x}^x)=S(G'/M^y,D_{M^y}^y)$ if
and only if there is a $g\in G'$ such that $M^{xg}=M^y$ and
$D_{M^{xg}}^{xg}\cong D_{M^y}^y$. This means that $x$ and $y$ are in the same 
$\cN(D_M)\times G'$-orbit of $\NMGd$, where
$\cN(D_M)=\stab_{\cN_G(M)}(D_M)$. Note that $\cN(D_M)$ acts indeed on 
$\NMGd$ by left translation since $\cN_G(M)$ does. On the
other hand, Theorem \ref{factor} tells us that every direct summand of
$\T^G_{G'}\bar{X}_G$ contains at most one composition factor from one 
HC-series $S(G'/M^x, D_M^x)$, and all such series occur, and we have
shown:

\begin{lemma}\mlabel{remcompfact} 
Keep the notation of Theorem \rmref{factor}. Let $x\in\cNNGdL$. Then 
$[x]_{G',X_L}\cap\NMGd$ is either empty or a union of 
$\cN(D_M)\times G'$-orbits in $\NMGd$. 
\end{lemma}

For general linear groups we have $\cN(D_M)\leq\cN(X_L)$ in the
relevant cases, and the result is trivially true. We do not know if
this holds in general. Theorem \ref{factor} implies now immediately:

\begin{cor}\mlabel{factorcor} Suppose $[x]_{G',X_L}\cap\NMGd$ is
  empty. Then no composition factor of $X_{G',x}$  lies in an
  HC-series $S(G'/M^y,D_{M^y})$ for some $y\in\NMGd$. If 
$[x]_{G',X_L}\cap\NMGd$ is the disjoint union of the
$\cN(D_M)\times(G'\cap N)$-orbits of $y_1,\ldots,y_m\in\NMGd$,
$$
[x]_{G',X_L}\cap\NMGd=\bigcup_{i=1}^m[y_i]_{G',D_M},
$$
then for every $1\leq i\leq m$ there is precisely one composition factor
$D_{G',M^{y_i}}$ of $X_{G',x}$ in the HC-series 
$S(G'/M^{y_i},D_M^{y_i})$.
\end{cor} 

In the special case $G=G'$ we have obviously
\begin{equation}\mlabel{factorspecial1}
\cnorm{G}{M}{G}(D_M)=\cnorm{G}{L}{G}=G,
\end{equation}
hence we have precisely one orbit
\begin{equation}\mlabel{factorspecial2}
[1]_{G,X_L} = [1]_{G,D_M} = G.
\end{equation}
This implies:

\begin{cor}\mlabel{factorspecial3} Let $M\in \L^c(X_G,Y_L)$. Then $X_G$
  has a unique composition factor $D_{G,M}$ in $S(G/M,D_M)$.
\end{cor}

The next theorem provides information how the composition factor $D_{G,M}$
behaves under HC-restriction:

\begin{thm}\mlabel{irredhcrest} Let $M\in \L^c(X_G,Y_L)$ and 
$M\leq G'\in\L_{G,L}$. Then
$$ 
\T^G_{G'}D_{G,M}=\bigoplus_{x\in\cnorm{G'}{M}{G}(D_M)}D_{G',M^x}.
$$ 
In particular, $\T^G_{G'}D_{G,M}$ is semisimple. Moreover, if $C$ is
any composition factor of $\bar{X}_G$ not
isomorphic to some $D_{G,M}$ then $\T^G_{G'}C$ has no composition factor equal
to some $D_{G',M^x}$ for $x\in\cnorm{G'}{M}{G}(D_M)$.
\end{thm}
\pf By definition we have $D_{G,M}\in S(G/M,D_M)_k$. Hence, by
Lemma \ref{dualstatement2}, we find irreducible $kG'$-modules $Z,Z'\in
S(G'/M,D_M)_k$ such that $Z$ is in the socle and $Z'$ is in the head of
$\T^G_{G'}D_{G,M}$. 

By Theorem \ref{factor}, $\T^G_{G'}\bar{X}_G$ has a unique composition 
factor in $S(G'/M,D_M)_k$, namely $D_{G',M}$. Now 
$D_{G,M}$ is a composition factor of $\bar{X}_G$ hence 
$\T^G_{G'}D_{G,M}$ is a subfactor of $\T^G_{G'}\bar{X}_G$. Thus every 
composition factor of $\T^G_{G'}D_{G,M}$ is a composition factor of 
$\T^G_{G'}\bar{X}_G$. We conclude that $\T^G_{G'}D_{G,M}$
has precisely one composition factor in $S(G'/M,D_M)_k$, namely 
$D_{G',M}$ and we have shown $Z=Z'=D_{G',M}$. Being in the head and 
in the socle of $\T^G_{G'}D_{G,M}$, the $kG'$-module $D_{G',M}$ splits
off. 

Applying Corollary \ref{factorcor} we get similarly that $D_{G',M^x}\in
S(G'/M^x, D_M^x)$ splits of $\T^G_{G'}D_{G,M}$ for every 
$x\in\cnorm{G'}{M}{G}(D_M)$. Thus 
$$
\T^G_{G'}D_{G,M}=\bigoplus_{x\in\cnorm{G'}{L}{G}(D_M)}D_{G',M^x}\oplus Y,
$$
where $Y$ is some $kG'$-module, whose composition factors do not
belong to any of  
the series $S(G'/M^x, D_M^x)$, $x\in N_{M\subseteq G'}$. But 
since $D_{G,M}\in S(G/M,D_M)$, every irreducible summand of the 
socle and the head of $\T^G_{G'}D_{G,M}$, hence of the socle and the 
head of $Y$, belongs to one of these HC-series. Therefore $Y=(0)$.

We see in particular that the composition factors of 
$\T^G_{G'}\bar{X}_G$ from these series are already composition factors
of the subfactor $\T^G_{G'}D_{G,M}$ of $\T^G_{G'}\bar{X}_G$ and the 
last assertion also follows.
\epf

\begin{rem}\mlabel{dirdecmod}
Corollary \rmref{factorcor} provides the
  recipe how the direct summands $D_{G',M^y}$ of $\T^G_{G'}D_{G,M}$ 
divide up into composition factors of the direct summands $X_{G',x}$
of $\T^G_{G'}X_G$, where $y$ runs through $\cnorm{G'}{M}{G}(D_M)$ and
$x$ through $\cNNML$.
\end{rem}

\begin{thm}\mlabel{homset} Let $M,G'\in \L_{G,L}$ such that 
$M\leq G'$. Then every
  homomorphism from $\R_M^{G'}Y_M$ to $X_{G'}$ factors through 
$\R_M^{G'}X_M$ and 
$$
\Hom_{\O G'}(\R_M^{G'}X_M, X_{G'})\simeq \O\simeq\Hom_{\O G'}
(X_{G'},\R_M^{G'}X_M).
$$
If $\phi=\phi_{G',M}$ in the first and $\rho=\rho_{G',M}$ in the second
$\Hom$-space are chosen such that every other homomorphism is a scalar
multiple, then the following holds, writing $1_k\otimes _{\O}\phi=
\bar{\phi}$ and $1_k\otimes _{\O}\rho=\bar{\rho}$:
\begin{enumerate}
\item The kernel of $\bar{\rho}$ has no composition factor isomorphic to
$D_{G',M'}$ for $M'\in \L^c(X_{G'},Y_L)$.
\item The image $U$ of $\bar{\phi}$ has simple head and contains 
$D_{G',M'}$ as composition factor for $M'\in \L^c(X_{G'},Y_L)$ if and 
only if $M'$ is conjugate in $G'$ to a subgroup of $M$.
\item If in addition $M\in \L^c(X_{G'},Y_{L})$ then $D_{G',M}$ is the head
$U/\Jac(U)$ of $U=\im(\bar{\phi})$. Moreover, if $V\leq \bar{X}_{G'}$ 
and $D_{G',M}$ is a
composition factor of $V$ and $M'\leq M,\,M'\in\L^c(X_{G'},Y_L)$, then
$D_{G',M'}$ is a composition factor of $V$. 
\end{enumerate}\end{thm}
\pf  In view of Result \ref{projrestsyssubgrp} we may assume that
$G=G'$. Now by Lemma \ref{functprop} every 
homomorphism from $\R_M^{G}Y_M$ to $X_{G}$ factors through
$\pi_{(L,\chi_L)}(\R_M^{G}Y_M)$. But if $\M'$ denotes the set of
$G$-conjugates of $L$ and $\C'$ the set of $G$-conjugate characters of
$\chi_L$, then $(\M',\C')$ is closed under conjugation in $G$. On the
other hand on $\O G$-lattices $\pi_{(\M',\C')} = \pi_{(L,\chi_L)}$ and
hence by Corollary \ref{commute} and Theorem \ref{projdec} 
$$
\pi_{(L,\chi_L)}(\R_M^GY_M) = \pi_{(\M',\C')}(\R_M^GY_M)
= \R_M^G\pi^M_{(\M',\C')}(Y_M) = \R_M^G(\pi_{(L,\chi_L)}(Y_M)) 
= \R_M^G(X_M).
$$
By Frobenius reciprocity and Theorem \ref{structhcrest}
\begin{equation}\mlabel{eq1}
\Hom_{\O G}(R_M^{G}X_M, X_{G})\simeq \bigoplus_{x\in\cnorm{M}{L}{G}}
\Hom_{\O M}(X_M,X_{M,x})
\simeq \Hom _{\O M}(X_M,X_M)\simeq \O,
\end{equation}
since $KX_{M,x}\not\simeq KX_M$ for $1\neq x\in\cnorm{M}{L}{G}$. 
Similarly
\begin{equation}\mlabel{eq2}
\Hom_{\O G}(X_{G},R_M^{G}X_M)\simeq \O.
\end{equation}
 We take now $\phi$ and $\rho$ to correspond to $1\in \O$ in Equation
 \ref{eq1} and \ref{eq2} respectively.  
We first prove part (ii). By Frobenius reciprocity we have
\begin{equation}\mlabel{eq3}
\begin{split}
k\simeq \Hom
_{kM}(\bar{Y}_M,\T^G_MU)&\simeq\Hom_{kG}(\R_M^G\bar{Y}_M,U)\\
&\cong\Hom_{kG}(\R_M^G\bar{Y}_M,\bar{X}_G)\\
&\simeq\Hom_{kG}(\R_M^G\bar{X}_M,\bar{X}_G)\\
&\simeq \Hom _{kM}(\bar{Y}_M,\bar{X}_M)
\end{split}\end{equation}
which shows in particular that $D_M$ is a composition factor of
$\T^G_MU$ and hence $\bar{X}_M$ is a submodule of $\T^G_MU$. Thus
$D_{M'}$ is a composition factor of $\T^G_{M'}U\ge\T^M_{M'}\bar{X}_M $
for every $M'\in\L^c(X_G,Y_L)$ with $M'\leq M$ and hence $D_{G,M'}$ is
a composition factor of $U$ by Theorem \ref{factor} (applied in the
special case $G=G'$). Since all these composition factors have
multiplicity one as composition factors of $\bar{X}_G$, the same holds
for the submodule $U$ of $\bar{X}_G$. Observe that the projective
$kG$-module $\R_M^G\bar{Y}_M$ contains the projective cover $Z$ of $U$ as
direct summand, and this is the projective cover of $\hd(U)$ as well. 
Since $\Hom_{kG}(\R_M^G\bar{Y}_M,U)$ is one dimensional, the same is
true for $\Hom_{kG}(Z,U)$ and $\hd(U)$ must be simple.

On the other hand let $M'\in \L^c(X_G,Y_L)$ and assume that $D_{G,M'}$
is a composition factor of $U$. Then it is a composition factor of 
$\R_M^G\bar{X}_M$ since $\bar{\phi}$ maps this module onto $U$.
Now $D_{G,M'}$ is in the HC-series $S(G/M',D_{M'})$ and therefore a
composition factor of $\hd(\R_{M'}^GD_{M'})$. Thus
$\R_{M'}^G\bar{Y}_{M'}$ contains the projective cover of
$D_{G,M'}$ and there is a nonzero homomorphism from
$\R_{M'}^G\bar{Y}_{M'}$ to $\R_M^G\bar{X}_M$. We have:

\begin{equation}\mlabel{eq4}\begin{split}
(0)&\neq\Hom_{kG}(\R_{M'}^G\bar{Y}_{M'},\R_M^G\bar{X}_M)\\
&\cong\Hom_{kM'}(\bar{Y}_{M'},\T^G_{M'}\R_M^G\bar{X}_M)\\
&\cong\bigoplus_{z\in\cD_{M,M}}\Hom_{kM'}(\bar{Y}_{M'},\R_{M'\cap
    M^z}^{M'}\T^{M^z}_{M'\cap M^z}(\bar{X}^z_{M^z})).
\end{split}\end{equation}
But the head of $\bar{Y}_{M'}$ is $D_{M'}$ and hence
cuspidal by assumption. Thus such a homomorphism can exist only if 
$M'\leq M^z$ for some $z\in G$ or equivalently $z^{-1}\in
N_{M'\subseteq M}$.  Applying Frobenius reciprocity to
\ref{eq4} again and conjugating we conclude that 
$$
\Hom _{kM^{\prime z^{-1}}}(\bar{Y}_{M^{\prime z^{-1}}}, \T^M_{M^{\prime
    z^{-1}}}\bar{X}_M)\neq(0).
$$
Setting $M^{\prime z^{-1}}=\tilde{M}$ we have by Corollary 
\ref{structhcrest}
\begin{equation*}\begin{split}
\Hom_{k\tilde{M}}(\bar{Y}_{\tilde{M}},\T^M_{\tilde{M}}\bar{X}_M)
&=\bigoplus\limits_{x\in\cnorm{\tilde{M}}{L}{M}}\Hom
_{k\tilde{M}}( \bar{Y}_{\tilde{M}},\bar{X}_{\tilde{M},x})\\
&=\Hom_{k\tilde{M}}( \bar{Y}_{\tilde{M}},\bar{X}_{\tilde{M}})=k.
\end{split}\end{equation*}
Thus we have $D_{G,\tilde{M}}\simeq D_{G,M'}$ since $\tilde{M}$ and
$M'$ are conjugate in $G$ and (ii) is shown.

To prove (iii) we observe that $D_{G,M}$ is a composition factor of $U$ by
part (ii). Suppose $D_{G,M}$ is not the head of $U$ then there is a proper
submodule $V$ of $U$ such that $D_{G,M}$ is a composition factor of $V$. As in
\ref{eq3} we show that $\Hom_{kM}(\bar{Y}_M,\T^G_MV)\neq(0)$, thus 
$\im\bar{\phi}\leq V$, a contradiction and therefore (iii) holds.

We now prove (i). First observe that 
$$
\R^G_M\bar{X}_M\simeq \bigoplus_x \bar{X}_M\otimes x,
$$
where $x$ runs through a set of coset repesentatives of a parabolic 
subgroup of $G$ having Levi complement $M$. Let 
$\tau:\R^G_M\bar{X}_M\rightarrow \bar{X}_M$ be the projection onto the 
summand $\bar{X}_M=\bar{X}_M\otimes 1$. By general
theory the map $\hat{\rho}$ defined to be the restriction of 
$\tau\circ\bar{\rho}$ to $T^G_M\bar{X}$
is $KM$-linear and is the map which corresponds to $\rho$ under the
isomorphism 
\begin{equation}\mlabel{eq5}
\Hom_{kG}(\bar{X}_G,\R^G_M\bar{X}_M)\simeq
\Hom_{kM}(\T^G_M\bar{X}_G,\bar{X}_M)
\end{equation}
In particular $\hat{\rho}$ is non-zero, indeed it is an epimorphism. By
Theorem \ref{structhcrest} and transitivity of HC-induction we have 
$$
\T^M_L\T^G_M \bar{X}_G=\T^G_L\bar{X}_G=\bigoplus\limits_{x\in
\cnorm{L}{L}{G}}\bar{X}_{L,x},
$$
in particular $\bar{X}_L$ occurs in $\T^G_L\bar{X}_G$ and hence in
$\T^M_L(\T^G_M\bar{X}_G)$ exactly once as a direct summand. Of course 
the same holds for $\bar{X}_M$, that is $\bar{X}_L$ occurs in 
$\T^M_L\bar{X}_M$ exactly once. Suppose that the unique subspace
$\bar{X}_L$ of $\bar{X}_G$ is in the kernel of $\bar{\rho}$ then it is also
in the kernel of $\tau\circ\bar{\rho}=\hat{\rho}$. But we may take
$\hat{\rho}$ to be the projection of 
$$
T^G_M\bar{X}_G=\bigoplus\limits_{x\in\cNNML}\bar{X}_{M,x}
$$
onto $\bar{X}_{M,1}\simeq \bar{X}_M$ and therefore the restriction
$\tilde{\rho}$ of 
$\hat{\rho}$ to $\bar{X}_{M,1}$ is the identity map. Thus 
$$
\T^M_L\tilde{\rho}:\T^M_L\bar{X}_{M,1}\rightarrow \bar{X}_{M}
$$
maps the unique subspace $\bar{X}_L$ of $\T^G_M\bar{X}_G$ into 
$\bar{X}_L\leq\T^M_L\bar{X}_M\leq\bar{X}_M$. Thus $\bar{X}_L$ can not 
be contained in the kernel of $\hat{\rho}=\tau\circ\bar{\rho}$ and 
hence not in the kernel of $\bar{\rho}$. Suppose $D_{G,M}$ is a
composition factor of $\ker(\bar{\rho})\leq X_G$. Then by (iii)
$D_{G,M'}$ is composition factor of $\ker(\bar{\rho})$ for every
$M'\in\L(X_G,Y_L)$ with $M'\leq M$. But $L\in\L^c(X_G,Y_L)$ by
\ref{irredmod}, hence $D_{G,L}$ is a composition factor of
$\ker(\bar{\rho})$. Thus we have a nonzero homomorphism from the
projective cover of $D_{G,L}$, which is a direct summand of
$\R_L^G\bar{Y}_L$ into $\ker(\bar{\rho})$, hence by Frobenius
reciprocity from $\bar{Y}_L$ into $\T^G_L\ker(\bar{\rho})\leq
\T^G_L\bar{X}_G$. But $\Hom_{kL}(\bar{Y}_L,\T^G_L\bar{X}_G)$ is one
dimensional, and its image is the direct summand $\bar{X}_L$ of
$\T^G_L\bar{X}_G$. We conclude that $\bar{X}_L$ is contained in
$\T^G_L\ker(\bar{\rho})$ and hence in $\ker(\bar{\rho})$, a
contradiction. This proves (i) and the theorem is shown.
\epf

\begin{cor}\mlabel{head} Let $G'\in \L_{G,L}$. Then the head 
$D_{G'}$ of $\bar{X}_{G'}$ is isomorphic to $D_{G',M}$ for some 
maximal $M\in \L^c(X_G,Y_L)$ with $M\leq G'$.
\end{cor}
\pf In view of Result \ref{projrestsyssubgrp} we may again assume that
$G=G'$. If $D_G$ is cuspidal we are done, since then
$G\in\L^c(X_G,Y_L)$ is maximal, and $D_G$ is the head of $\bar{X}_G$. 

Assume $G\notin\L^c(X_G,Y_L)$ and let $M\in \L_G$ be the semisimple
HC-vertex of $D_G$. We show that $M\in \L_{G,L}$. By general
assumption \ref{defrestsys} $\chi_{G}$ occurs as constituent of 
$\R^{G}_L\chi_L$. Hence $\T^G_M\chi_{G'}$ is a summand of 
$$
\T^G_M\R^G_L\chi_L=\sum_{x\in\D_{L,M}}\R^M_{L^x\cap M}\T^{L^x}_{L^x\cap
  M}\chi^x_{L^x}.
$$
But $\T^{L^x}_{L^x\cap M}\chi^x_{L^x}=(0)$ unless $L^x\leq M$. This 
shows that for every composition factor $S$ of $\bar{X}_G$ we have 
$\T^G_MS=(0)$ unless $L^x\leq M$ for
some $x\in G$. As $X_G$ and $X^x_G$ are isomorphic for
such an $x$, we may assume that the semisimple HC-vertex of $D_G$ 
contains $L$. 

The short exact sequence 
$$
0\to \Jac(\bar{X}_G)\to \bar{X}_G\to D_G\to 0
$$
yields a short exact sequence 
$$
0\to\T^G_M\Jac(\bar{X}_G)\to \bigoplus_{y\in\cnorm{M}{L}{G}}
\bar{X}_{M,y}\to\T^G_MD_G\to 0
$$
using Theorem \ref{structhcrest}. Consequently an irreducible direct
summand $S$
of the head of $\T^G_MD_G$ has to be the head of one
of the summands of the middle term, and thus by induction it has to be of the
form $D_{M,y}$ for some $y\in\cnorm{M}{L}{G}$. By Lemma 
\ref{dualstatement2}, $D_{M,y}$ is cuspidal and $D_G$ is in the
HC-series $S(G/M,D_{M,y})_k$.
But $D_{M,y}=(D_{M^{y^{-1}}})^y$ by \ref{2.9a}. Thus replacing $M$ by
$M^{y{-1}}$ if necessary we may assume $y=1$. By Theorem
\ref{factor} the unique composition factor of $\bar{X}_G$ in
the series $S(G/M, D_M)$ is $D_{G,M}$. Theorem \ref{homset} part
(iii) tells us now that $M$ has to be maximal with $M\in \L^c(X_G,Y_L)$.
\epf

Recall our definitions in \ref{tauandt}. We apply Theorem \ref{homset}
in the special situation $M=L$: 

\begin{cor}\mlabel{sublat} Let $G'\in \L_{G,L}$. Then 
$$
\rho_{G'} =\rho_{G',L}:X_{G'}\to\R_L^{G'}X_L
$$
is injective and its image $X_\rho = X_{\rho,G'} = \rho_{G'}(X_{G'})$
is a sublattice of $X'_{G'}\leq\R_L^{G'}X_L $ satisfying $KX_\rho =
KX'_{G'}$ and $X'_{G'}=\sqrt{\rho(X_{G'})}$. Morover every
homomorphism from $Y_{G'}$ into
$\R_L^{G'}X_L$ factors through a multiple of $\rho$ and is
the combined image 
$$
X_\rho=\tau_{Y_{G'}}(\R_L^{G'}X_L)=\tau_{Y_{G'}}(X'_{G'})
$$
(which means  $X_\rho$ is the combined image 
of $Y_{G'}$ in $\R_L^{G'}X_L$ and $X_{G'}$ respectively). Finally the
composition factors of $\ker\bar{\rho}$ and the cokernel 
$V=\coker\bar{\rho}=\bar{X}'_{G'}/\bar{X}_\rho$ coincide
(multiplicities counted). In particular no composition factor of $V$
has the form $D_{G',M'}$ for some $M'\in\L^c(X_{G'},Y_L)$.
\end{cor}     
\pf By Theorem \ref{projdec} and Frobenius reciprocity 
\begin{equation*}\begin{split}
\Hom _{\O G'}(Y_{G'},\R_L^{G'}X_L)&\simeq \bigoplus \limits_{x\in
\cnorm{L}{L}{G}}\Hom_{\O L}(Y_{L,x},X_L)\oplus \Hom_{\O L}(Z,X_L)\\
&\simeq \Hom_{\O L}(Y_{L},X_L)\\
&\simeq \O,
\end{split}\end{equation*}
since $X_{L,x}\not\simeq X_L$ for $1\not=x\in\cnorm{L}{L}{G}$ by
Theorem \ref{structhcrest}, $\pi_{(L,\chi_L)}(Y_L)=X_L$ and since
$X_L$ is cuspidal, by our basic assumption \ref{notrestsys}.  Note
that $ \pi_{(L,\chi_L)}(Z)=(0)$ by Theorem \ref{projdec}. Thus there is a
unique homomorphism $\hat{\rho}:Y_{G'}\to  R_L^{G'}X_L$ corresponding
to $1\in \O$. Again by Lemma \ref{functprop}
$$
\Hom _{\O G'}(Y_{G'},\R_L^{G'}X_L)\simeq\Hom _{\O
  G'}(X_{G'},\R_L^{G'}X_L),
$$
thus every homomorphism from $Y_{G'}$ to $\R_L^{G'}X_L $ factors through
$X_{G'}$, and we can choose $\rho:X_{G'}\to\R_L^{G'}X_L $ such that
$\hat{\rho}$ is the composite of the natural projection of $Y_{G'}$
onto $X_{G'}=\pi^{G'}(Y_{G'})$ and  $\rho$. Thus $\im\rho$ is the
combined image of $Y_{G'}$ in $\R_L^{G'}X_L$ and in every submodule of
$\R_L^{G'}X_L$ containing $\im\rho$. Now
$1_K\otimes _{\O} \rho:KX_{G'}\to\R_L^{G'}KX_L $ has to be injective
since $KX_{G'}$ is irreducible by Theorem \ref{structhcrest}, hence $\rho$ is
injective, too. Thus $1_K\otimes_\O\rho$ maps $KX_{G'}$ isomorphically
onto $KX'_{G'}\leq K\R_L^{G'}X_L$.  In particular, since 
$$
X'_{G'} = \sqrt{X'_{G'}} =
KX'_{G'}\cap\R_L^{G'}X_L=K\rho(X_{G'})\cap\R_L^{G'}X_L = \sqrt{\rho(X_{G'})}
$$
we have 
$$
\rho(X_{G'})\leq X'_{G'}.
$$
Since $KX_{G'}=KX'_{G'}$ the composition factors of $\bar{X}_{G'}$ and
$\bar{X}'_{G'}$ coincide. Now by the first isomorphism theorem 
$$
\bar{X}_{G'}/\ker\bar{\rho}\cong\im\bar{\rho}\leq\bar{X}'_{G'}.
$$
Consequently $\ker\bar{\rho}$ and $\coker\bar{\rho}$ have the same
composition factors. 
\epf

\begin{lemma}\mlabel{compfact} Let $G'\in \L_{G,L}$ and let $S$ be a
  composition factor of $\bar{X}_{G'}$ with $S\not =D_{G',M'}$ for all
  $M'\in \L^c(X_{G'},Y_L)$. Let $M\in\L_{G,L}$. Then
$$
\Hom_{k\G}(\R_M^\G\bar{Y}_M,\R_{G'}^\G S) = (0).
$$
Thus no composition factor of $\R_{G'}^\G S$ is contained in
the HC-series $S(\G/M,D_M)_k$. 
\end{lemma}
\pf 
We have
\begin{equation}\mlabel{compfacteq3}\begin{split}
\Hom_{k\G}(\R_M^\G\bar{Y}_M,\R_{G'}^\G S) 
&=\Hom_{kM}(\bar{Y}_M,\T^G_M\R_{G'}^GS)\\
&=\bigoplus_{z\in \tilde{\cD}_{MG'}\cap \tilNLM}\Hom_{kM}(\bar{Y}_M,
\R_{M\cap {G'}^z}^M\T^{{G'}^z}_{M\cap {G'}^z}S^z)\\
&=\bigoplus_{z\in \tilde{\cD}_{MG'}\cap \tilNLM}\Hom_{k(M\cap {G'}^z)}
(\T^M_{M\cap G^{\prime z}}\bar{Y}_M,\T^{G^{\prime z}}_{M\cap G^{\prime
    z}}S^z),
\end{split}\end{equation}
where $\tilde{\cD}_{MG'}$ denotes a suitable system of double coset
representatives of parabolic subgroups of $\G$ containing $G'$ and $M$
as Levi complements. Now
\begin{multline}\label{compfacteq4}
\Hom_{k(M\cap G^{\prime z})}
(\T^M_{M\cap G^{\prime z}}\bar{Y}_M,\T^{G^{\prime z}}_{M\cap G^{\prime
    z}}S^z)\\
=\bigoplus_{x\in\cnorm{M\cap G^{\prime z}}{L^z}{M}} 
\Hom_{k(M\cap G^{\prime z})}(\bar{Y}_{M\cap G^{\prime z},x},
\T^{G^{\prime z}}_{M\cap G^{\prime z}}S^z)\oplus
\Hom_{k(M\cap G^{\prime z})}(\bar{Z}_z,
\T^{G^{\prime z}}_{M\cap G^{\prime z}}S^z),
\end{multline}
where $Z_z$ is a projective $\Or M\cap G^{\prime z}$-module with 
$\pi^{M\cap G^{\prime z}}(Z_z)=(0)$ by Theorem \ref{projdec}.

Now $S^z$ is a subfactor of $\bar{X}_{G'}^z$. By Corollary
\ref{hcrestx}, 
applied to the projective restriction system 
$\PR^z(X_G,Y_L)=\PR(X^z_G,Y^z_L)$, we have
$$
\T^{G^{\prime z}}_{M\cap G^{\prime z}}X_{G'}^z=
\bigoplus_{y\in\cnorm{M\cap G^{\prime z}}{L^z}{G^{\prime z}}} 
X^{(z)}_{M\cap G^{\prime z},y},
$$
with 
$X^{(z)}_{M\cap G^{\prime z},y}=
\pi_{(L^{zy},\chi^{zy})}
(\T^{G^{\prime z}}_{M\cap G^{\prime z}}X_{G'}^z)$ 
(see Notation \ref{moreprojnot}). Similarly we denote the projective cover of 
$X^{(z)}_{M\cap G^{\prime z},y}$ by $Y^{(z)}_{M\cap G^{\prime z},y}$
and we 
write $D^{(z)}_{M\cap G^{\prime z},M^{\prime x}}$ for the composition
factors 
of $\T^{G^{\prime z}}_{M\cap G^{\prime z}}X_{G'}^z$
lying in series $S((M\cap G^{\prime z})/M^{\prime x},D_{M^{\prime
    x}})$ for 
$M^{\prime x}\in \cL^c(X_{G^{\prime}}^z,Y_L^z)$.

From Theorem \ref{projdec} and Corollary \ref{commute} follows 
$$
\T^{G^{\prime z}}_{M\cap G^{\prime z}}\bar{X}_{G'}^z =
\T^{G^{\prime z}}_{M\cap G^{\prime z}}\pi^{M\cap G^{\prime
    z}}\bar{Y}_{G'}^z=
\pi^{G^{\prime z}}\T^{G^{\prime z}}_{M\cap G^{\prime z}}\bar{Y}_{G'}^z,
$$
and hence by Lemma \ref{functpropproj}  
\begin{equation*}\begin{split}
\Hom_{k(M\cap G^{\prime z})}
(\bar{Z}_z,\T^{G^{\prime z}}_{M\cap G^{\prime z}}S^z)&\le 
\Hom_{k(M\cap G^{\prime z})}
(\bar{Z}_z,\T^{G^{\prime z}}_{M\cap G^{\prime z}}\bar{X}_{G'}^z)\\
&\cong \Hom_{k(M\cap G^{\prime z})}(\bar{Z}_z,
\pi^{G^{\prime z}}\T^{G^{\prime z}}_{M\cap G^{\prime
    z}}\bar{Y}_{G'}^z)\\
&\cong \Hom_{k(M\cap G^{\prime z})}(\pi^{G^{\prime z}}\bar{Z}_z,
\T^{G^{\prime z}}_{M\cap G^{\prime z}}\bar{X}_{G'}^z)\\
&\cong \Hom_{k(M\cap G^{\prime z})}((0),
\T^{G^{\prime z}}_{M\cap G^{\prime z}}\bar{X}_{G'}^z)\\
&=(0).
\end{split}\end{equation*}

Now suppose that
$$
\Hom_{k(M\cap G^{\prime z})}
(\bar{Y}_{M\cap G^{\prime z},x},\T^{G^{\prime z}}_{M\cap G^{\prime z}}S^z)\le 
\bigoplus_{y\in\cnorm{M\cap G^{\prime z}}{L^z}{G^{\prime z}}} 
\Hom_{k(M\cap G^{\prime z})}(\bar{Y}_{M\cap G^{\prime z},x},
\bar{X}^{(z)}_{M\cap G^{\prime z},y})\not=(0),
$$
say $\Hom_{k(M\cap G^{\prime z})}(\bar{Y}_{M\cap G^{\prime z},x},
\bar{X}^{(z)}_{M\cap G^{\prime z},y_0})\not=(0)$. 
Then using again Theorem \ref{projdec} we see that 
$\bar{Y}_{M\cap G^{\prime z},x}$ is the projective cover of 
$\bar{X}^{(z)}_{M\cap G^{\prime z},y_0}$. 
That is, we have $\bar{Y}_{M\cap G^{\prime z},x}\cong 
\bar{Y}^{(z)}_{M\cap G^{\prime z},y_0}$,
hence its head is one of the modules 
$D^{(z)}_{M\cap G^{\prime z},M^{\prime x}}$ for some 
$M^{\prime x}\in \cL^c_{G^{\prime z},L^z}$ by Corollary \ref{head}. 
Now Theorem \ref{irredhcrest} implies 
$\Hom_{k(M\cap G^{\prime z})}(\bar{Y}_{M\cap G^{\prime z},x},
\T^{G^{\prime z}}_{M\cap G^{\prime z}}S^z)=(0)$ and the lemma follows. 
\epf

In the next result we extend Corollary \ref{sublat} to HC-induced
lattices. We define
\begin{equation}
\bY=\bY^{\G}=\bigoplus_{G'\in  \L_{{G},L}}\R_{G'}^\G Y_{G'}.
\end{equation}

\begin{cor}\label{byspur} Let $G'\in \L_{G,L}$ and let 
$X_\rho\leq X'_{G'}\leq \R_L^{G'}$ be defined as in Corollary
\rmref{sublat}. Then
$$
\R_{G'}^\G X_\rho = \tau_\bY(\R_{G'}^\G X'_{G'}). 
$$
\end{cor}

\pf By \ref{sublat} $X_\rho$ is image of $Y_{G'}$ hence $\R_{G'}^\G
X_\rho$ is epimorphic image of $\R_{G'}^\G Y_{G'}$, which is a direct
summand of $\bY$. We conclude that  
$$
\R_{G'}^\G X_\rho \leq  \tau_\bY(\R_{G'}^\G X'_{G'}). 
$$
Let $S$ be a composition factor of the cokernel of $\bar{\rho}$. Then,
by Lemma \ref{sublat} $S$ is not of the form $D_{G',M'}$ for some 
$M'\in\L^c(X_{G'},Y_L)$. Lemma \ref{compfact} tells us
$$
\Hom_{k\G}(\R_M^\G\bar{Y}_M,\R_{G'}^\G S) = (0).
$$
Thus the image of every homorphism from $\bar{\bY}$ to
$\R_{G'}^\G\bar{X}_{G'}$ has to be contained in $\im\bar{\rho}$. We
lift this to $\O$ to get the desired result.
\epf

\begin{lemma}\mlabel{lattiso} 
Suppose that every composition factor of $\R_L^G\bar{X}_L$
  is isomorphic to some $D_{G,M}$ where $M\in \L^c(X_G,Y_L)$. Then we have
\begin{enumerate}
\item The $kL$-module $\bar{X}_L$ is irreducible, that is
  $\bar{X}_L=D_L$.\\
\item Let $G'\in \L_{G,L}$. Then every composition factor of
  $R_L^{G'}\bar{X}_L$ is of the form $D_{G',M}$ for some
  $\L^c(X_{G'},Y_{L})=\{M\in\L^c(X_G,Y_L)\mid M\leq G'\}$.\\
\item There is an isomorphism between the poset of local submodules of
$\bar{X}_{G'}$ and the poset of subgroups $\L^c(X_{G'},Y_{L})$ of $G'$.\\
\item $R_{G'}^{\tG}X_{G'}$ is a pure sublattice of $R_L^{\tG}\chi_L$
isomorphic to $R_{G'}^{\tG}X_{G'}'$ for all $G'\in \L_{G,L}$.
\end{enumerate}\end{lemma}

\pf Lemma \ref{compfact} implies immediately part (i) and Theorem
\ref{homset} part (iii) whence we have established part (ii) observing that
all composition factors of $\bar{X}$ are composition factors of
$R_L^G\bar{X}_L$, since $KX_G$ is a constituent of $R_L^{G'}X_L$ by
Definition \ref{defrestsys}.

To show part (ii) let $D$ be a composition factor of
$R_L^{G'}\bar{X}_L$ and let $D\in S(G'/M',E)_k$, where $E$ is an
irreducible cuspidal $kM'$-module and $M'\in \L_{G'}$. Since $D$ is a 
composition factor of $R_L^{G'}\bar{X}_L$ and $X\in\lat_{\O L}$ is
cuspidal, we may choose $M'$ such that $L\leq M'$, that is $M'\in
\L_{G',L}$. Now $R_{G'}^G(D)$ is a subfactor of $R_L^G(\bar{X}_L)$
and, by Lemma \ref{dualstatement1}, every composition factor of the
head of $R_{G'}^G(D)$ is in $S(G/M',E)_k$. Our general assumptions
imply now that $S(G/M',E)_k=S(G/M,D_M)_k$ for some
$M\in\L^c(X_G,Y_L)$. Thus there exists $x\in N$ such that $M'=M^x$ and
$D_M^x=E$. Since $M'\in \L_{G,L}$ we have $M'\in \L^c(X_G,Y_L)$ and 
$E=D_{M'}$. But then $D\in S(G'/M', D_{M'})$ as desired.

It remains to show (iv): In Theorem \ref{homset} we constructed a
homomorphism $0\not= \rho_{G',L}:X_{G'}\to R_L^{G'}X_L$, which is
injective since $KG'$ is an irreducible $KG'$-module. Thus
$X_{G'}\subset 
R_L^{G'}X_L$. But part $ii)$ above and part $i)$ of Theorem \ref{homset} imply
that $\bar{\rho}_{G',L}:=1_k\otimes_{\O}\rho_{G',L}:\bar{X}_{G'}\to
R_L^{G'}\bar{X}_L$ is injective, too, hence $X_{G'}$ is a pure sublattice of 
$R_L^{G'}X_L$. By construction $X_{G'}'$ is a pure sublattice of $R_L^{G'}X_L$
affording the same character $\chi_{G'}$. Since the multiplicity of
$\chi_{G'}$ in $R_L^{G'}\chi_L$ is one, $X_{G'}'$ has to be the image of
$\rho_{G',L}$. Since HC-induction preserves purity and isomorphisms
part (iv) follows \epf\\

We point out that the assumption in Corollary \ref{lattiso} is frequently
satisfied, as for example in the application for $GL_n(q)$ and we shall see in
the following section how the conclusion of Corollary \ref{lattiso} can be
used to get information on decomposition numbers of $\tG$. 
\section{Hecke- and Schur algebras}\mlabel{schur}
We continue with the set up of the previous section. So let $\PR=\PR(X_G,Y_L)$
be a projective restriction system, consisting of the data
$\{X_M,Y_M,X_M'| M\in \L_{G,L}\}$. We define $Q=Q_L=\R_L^\G Y_L$ and 
$\be^\G=\be=\R_L^{\G}(\psi)$, where $\psi:Y_L\to\pi^L(Y_L)=X_L$ is the
natural projection.

\begin{thm}\mlabel{decmat}
The epimorphism $\be:Q\to\R_L^{\G}X_L$ satisfies
Hypothesis  \rmref{hypo}. Thus in particular the endomorphism ring
of $\R_L^{\G}\bar{X}_L$ is liftable and the decomposition matrix of
the endomorphism ring $\End_{\O\G}(\R_L^{\G}X_L)$ is part of
the $\ell$-modular decomposition matrix of $\G$.
\end{thm}

\pf By Corollary \ref{commute} and by \ref{notrestsys}
$$
\pi^\G(Q) = \pi^\G(\R_L^\G Y_L)=\R_L^\G\pi^L(Y_L) = \R_L^\G X_L,
$$
hence the theorem follows from Result \ref{projhom}.
\epf

Lemma \ref{ptl} implies now immediately:

\begin{cor}\mlabel{torsionless} The $\O\G$-lattice $\R_L^{G'}X_L$ is 
$Q$-torsionless.
\end{cor}
%

\begin{notation}\mlabel{assochecke}{\rm
We call $\H_R^{\G}=\H^{\G}_R(\PR)=\End_{R\G}(R\otimes
_{\O}\R_L^{\G}X_L)$ the \emph{Hecke algebra} over $R$ \emph{associated} with 
$\PR=\PR(X_{G},Y_L)$ on level $\G$. If no ambiguities arise we omit
super- and subscripts, for instance $\H_R$ is until further notice 
$\H_R^{\G}$. Observe that 
Theorem \ref{decmat} implies that $\H_R = R\H_\O=R\otimes_\O\H_\O$.
According to the general theory of quotients of Hom-functors as
outlined in the first section, we have associated with $\H_R$
quotients of Hom-functors $H_R$ given by
$$
H_R:\mod_{R\G}\to\mod_{\H_R}:V\to\Hom_{R\G}(Q,V)/\Hom_{R\G}(Q,V)J_\be,
$$
where $J_\be$ is the ideal of $\End_{R\G}(Q)$ consisting of those
endomorphisms of $Q$, whose image is contained in the kernel of $\be$}.
\end{notation}


Let $G'\in \L_{G,L}$. We have seen in Lemma \ref{sublat} that the irreducible
character $\chi_{G'}$ of $X_{G'}$ occurs exactly with multiplicity one in
$\R_L^{G'}\chi_L$. By Fittings Lemma we find a one dimensional
$\H_K^{G'}$-module generated by $\sigma =\sigma_{G'}$ such that
$K\sigma \R_L^{G'}X_L=KX_{G'}$. We choose $\sigma \in \H_\O^{G'}$ such
that $\sigma\H_\O^{G'}$ is pure in $\H_\O^{G'}$ and is
a generator of this one dimensional $\O$-space. 
Then $\sigma \R_L^{G'}X_L$
is an $\O$-lattice in $KX_{G'}$, but in general 
$\sigma \R_L^{G'}X_L\not\cong X_{G'}$ and $\sigma
\R_L^{G'}X_L\not\cong
X_{G'}'\leq \R_L^{G'}X_L$. 
However $X_{G'}'=\iota_{(G',\chi_{G'})}(\R_L^{G'}X_L)$ is pure in 
$\R_L^{G'}X_L$ by construction and affords $\chi_{G'}$. This implies:

\begin{lemma}\label{primepure} Let $U\le\R_L^{G'}X_L$ with $KU=KX'_{G'}$. Then
$$
U\le X'_{G'}=\sqrt{U}.
$$
In particular 
$$
\sigma \R_L^{G'}X_L\le
X_{G'}'\quad\text{and}\quad  \sqrt{\sigma \R_L^{G'}X_L}=X_{G'}'.
$$
\end{lemma}

\pf We have 
$$
\sqrt{U}=KU\cap\R_L^GX_L=KX'_G\cap\R_L^GX_L=X_G'.
$$ 
\epf

The next result relates $\H_{\O}^{G'}$ and $\H_{\O}^{\G}$: 

\begin{thm}\mlabel{relate}
Let $G'\in \L_{\G,L}$. Then the following
  holds for $R\in \{K,\O,k\}$:
\begin{enumerate}
\item $\H_R^{G'}$ is a subalgebra of $\H_R^{\G}$ and 
$\H_R^{\G}$ is
free as left $\H_R^{G'}$-module.
\item If I is a right ideal of $\H_{\O}^{G'}$ and $Z=I\R_L^{G'}X_L$ then 
$$
\R_{G'}^{\G}Z=I\H_{\O}^{\G}\R_L^{\G}
X_L=I\R_L^{\G}X_L,
$$
where $I\H_{\O}^{\G}\cong
I\otimes_{\H_{\O}^{G'}}\H_{\O}^{\G}$. In particular, if I is pure in
$\H_{\O}^{G'}$ then $I\H_{\O}^{\G}$ is pure in
$\H_{\O}^{\G}$.\end{enumerate}
In particular, the lattice $\sigma_{G'}\H_\O^\G$ is a pure right ideal
of $\H_O^\G$.
\end{thm}

\pf By Frobenius reciprocity
$$
\H_{\O}^{\G}=\Hom_{\O\G}
(\R_L^{\G}X_L,\R_L^{\G}X_L)=
\Hom_{\O G'}(\R_L^{{G'}}X_L,\T_{G'}^{\G}\R_L^{\G}X_L).
$$
Note that maps in the right hand side Hom-set correspond to maps in the left 
{\rm Hom}-set simply by extending the map in the natural way to
$\R_L^{\G}X_L=\R_L^{G'}X_L\O \G$, where $\R_L^{G'}X_L$ is
considered as module for a parabolic subgroup whose Levi complement is
$G'$. We choose a set $\tD_{L,G'}$ of double coset representatives of
parabolic subgroups of $\G$ having $L$ and $G'$ as Levi complements. 
Now by Mackey decomposition (see \ref{mackey}) 
\begin{equation}\begin{split}
 \Hom_{\O G'}(\R_L^{{G'}}X_L,\T_{G'}^{\G}\R_L^{\G}X_L)&\cong
\bigoplus_{x\in\tD_{L,G'}}\Hom_{\O G'}(\R_L^{{G'}}X_L,\R_{G'\cap
  L^x}^{G'}\T^{L^x}_{G'\cap L^x}X_L^x)\\
&\cong \bigoplus_{x\in\tD_{L,G'}\cap\NLGd}\Hom_{\O
  G'}(\R_L^{{G'}}X_L,\R_{L^x}^{G'}X_L^x),
\end{split}\end{equation}
since $X_L$ is cuspidal. But for $x=1$ we get $\H_{\O}^{G'\prime }=\Hom_{\O
  G'}(\R_L^{{G'}}X_L,\R_L^{{G'}}X_L)$ and for $1\not=x\in\tD_{L,G'}$, $L^x\le
  G'$, $\Hom_{\O G'}(\R_L^{{G'}}X_L,\R_{L^x}^{{G'}}X_L^x)$ is either
  $(0)$, if $L^x\not=_{G'} L$, or $L^{xy}=L$ for a $y\in G'$ but 
$X_L^{xy}\not\cong X_L$. Or it is as left $\H_{\O}^{G'}$-module isomorphic
to the regular module. A basis of this is then given for instance
by choosing one isomorphism from $\R_L^{G'}X_L$ to $\R_{L^x}^{G'}X_L^x$ 
induced by conjugation certain $y\in G'$. For $R=k, K$ we argue
analogously. Part ii) follows now easily.\epf

\begin{notation}\mlabel{notationlattices}{\rm 
Recall the definition of $\rho_{G'}:X_{G'}\to\R_L^{G'}X_L$ in
Corollary \ref{sublat}. 
The following $\O\G$-lattices will be considered:
\begin{alignat*}{3}
&\bX'&=&\bX^{\prime\G}&=&\bigoplus_{G'\in\L_{{G},L}}
\R_{G'}^{\G}X_{G'}'\\
&\bX&=&\bX^{\G}&=&\bigoplus_{G'\in\L_{{G},L}}\R_{G'}^{\G}X_{G'}\\
&\bX_h&=&\bX_h^{\G}&=&\bigoplus _{G'\in
  \L_{{G},L}}\sigma_{G'}\R_{G'}^{\G}X_L\\
&\bY&=&\bY^{\G}&=&\bigoplus_{G'\in  \L_{{G},L}}Q_{G'},
\end{alignat*}
setting $Q^\G_{G'}=Q_{G'}= \R_{G'}^\G Y_{G'}$.
Define $\brho:{\bf X}\rightarrow {\bf X}'$ to be the homomorphism
obtained by summing up the maps $\R_{G'}^\G\rho_{G'}$, for
$G'\in\L_{G,L}$. As in Corollary \ref{sublat} we denote the image of
$\rho=\rho_{G'}$ in $X'_{G'}$ by $X_\rho=X_{\rho_{G'}}$. Note that
$kX_\rho=\bar{X}_\rho$ is not the image of
$\bar{\rho}=1_k\otimes_\O\rho$ in general, but the latter is always an
epimorhic image of the former.} 
\end{notation}
 
Here is our main result:

\begin{thm}\label{bhomfunctor} Let the $\O\G$-lattices $\bX=\bX^\G$
  and $\bY=\bY^\G$ be defined as in \rmref{notationlattices}. Then
$$
\pi^\G(\bY)= \bX.
$$
Let $\bbe:\bY\to\bX$ denote the corresponding epimorphism. Then $\bbe$
satisfies Hypothesis \rmref{hypo} and the requirements of Theorem
\rmref{decmatrix}. In particular, the decomposition matrix of
$\End_{\O\G}(\bX)$ is part of the $\ell$-modular decomposition matrix
of $\G$.
\end{thm} 

\pf Corollary \ref{commut} implies
$$
\pi^\G(Q_M) = \pi^\G(\R_M^\G Y_M) = \R_M^\G X_M
$$
for every $M\in\L_{G,L}$. Since $\pi^\G$ obviously preserves direct
sums the theorem follows using Corollary \ref{projhom}.
\epf   

Lemma \ref{primepure} in conjunction with Theorem \ref{hom1corr} imply
that $\End_{\O\G}(\bX_h)=\End_{\O\G}(\bX')$ (see \ref{puresub}
below). Now Lemma \ref{lattiso} implies immediately:

\begin{thm}\label{lattiso2} Suppose that every composition factor of 
$\R_L^G\bar{X}_L$ is isomorphic to some $D_{G,M}$ where 
$M\in \L^c(X_G,Y_L)$. Then $\bX\cong\bX_\brho=\bX'$.
\end{thm}

So under the assumption of the theorem we can immediately compute the
decoposition of irreducible lattices occurring in $\bX$ in terms of
the Hecke algebra $\H_\O$.

Obviously we would now like to know $\End_{\O\G}(\bX)$ in general. We
shall use the following result to show, that we can relate this
endomorphism ring to the
endomorphism rings of the other lattices defined in \ref{notationlattices}.

\begin{thm}\mlabel{sublattices} Keep the notation introduced above.
Let $G'\in \L_{{G},L}$.\begin{enumerate}\item
We have $\sigma_{G'}\R_L^\G X_L\le\R_{G'}^\G X_\rho\le\R_{G'}^\G
X_{G'}'$, and therefore
$$
\bX_h\le \bX_\brho\le \bX',
$$
where $\bX_\brho=\bX_\brho^\G=\brho(\bX)$.
\item $\sigma_{G'}\R_L^\G X_L= \tau_{Q_L}(\R_{G'}^\G X_\rho)
=\tau_{Q_L}(\R_{G'}^\G X_{G'}')$ and
$\tau_{Q_{G'}}(\R_{G'}^\G X_{G'})=\R_{G'}^\G X_\rho$,
hence
$$
\tau_{Q_L}(\bX')=\bX_h=\tau_{Q_L}(\bX_\brho)\quad\text{and}\quad 
\tau_\bY(\bX') =\bX_\brho.
$$
\end{enumerate}\end{thm}

\pf Note that $K\sigma_{G'}\R_L^{G'}X_L=KX_\rho=KX'_{G'}$. Thus Lemma
\ref{primepure} implies 
$$
\sigma_{G'}\R_L^\G X_L\le\R_{G'}^\G X_{G'}\quad\text{and}\quad 
\R_{G'}^\G X_\rho\le\R_{G'}^\G X_{G'}',
$$ 
since HC-induction preserves injections. 

Corollary \ref{byspur} implies immediately
\begin{equation}\label{qmspur}
\tau_{Q_M}(\R_{G'}^\G X'_{G'})\leq\R_{G'}^\G X_\rho=
\tau_\bY(\R_{G'}^\G X'_{G'})
\end{equation}
for every  $M\in\L_{G,L}$, since $Q_M$ is a direct summand of $\bY$. 
Taking $M=L$ we obtain
$$
\tau_{Q_L}(\R_{G'}^\G X'_{G'})\leq\R_{G'}^\G X_\rho.
$$
Now $\sigma_{G'}\R_L^{G'}X_L$ is the image of the endomorphism 
$\sigma$ of $\R_L^{G'}X_L$ and hence the image of the surjective
composite map 
$$
\R_L^{G'}Y_L\to\R_L^{G'}X_L\stackrel{\sigma_{G'}}{\to}
\sigma_{G'}\R_L^{G'}X_L, 
$$
hence $\sigma_{G'}\R_L^\G X_L$ is the image of the composite map
$$
Q_L=\R_L^\G Y_L\to\R_L^\G X_L\stackrel{\hat{\sigma}}{\to}
\sigma_{G'}\R_L^\G X_L
$$ 
setting $\hat{\sigma}=\R_{G'}^\G(\sigma_{G'})$. We conclude that 
$$
\sigma_{G'}\R_L^\G X_L\leq\tau_{Q_L}(\R_{G'}^\G X'_{G'})\leq\R_{G'}^\G
X_\rho
$$
and
$$
\sigma_{G'}\R_L^\G X_L=\tau_Q(\sigma_{G'}\R_L^\G X_L).
$$
Summing up we get part i). 
By Theorem \ref{relate} $\sigma_{G'}\H_\O^\G$  is pure in $\H_\O^\G$, hence
by \cite[4.13]{di} and Lemma \ref{primepure} we have
$$
\sigma_{G'}\R_L^\G X_L=\tau_{Q_L}(\sigma_{G'}\R_L^\G X_L)=
\tau_{Q_L}(\sqrt{\sigma_{G'}\R_L^\G X_L})=\tau_{Q_L}(\R_{G'}^\G X'_{G'}).
$$
Summing up we get
$$
\tau_{Q_L}(\bX')=\bX_h=\tau_{Q_L}(\bX_\brho).
$$
Applying formula \ref{qmspur} for $M=G'$ we obtain
$$
\tau_{Q_{G'}}(\R_{G'}^\G X_{G'})\leq\R_{G'}^\G X_\rho
$$ 
but obviously we have equality here, since $\R_{G'}^\G X_\rho$ is
epimorphic image of $Q_{G'}=\R_{G'}^\G Y_{G'}$. Summing up we get
$$
\tau_\bY(\bX')=\bX_\brho.
$$
\epf

The additional hypothesis on $\H^\G_\O$ (see \ref{isym}) in the
following result is satisfied in all our applications:

\begin{thm}\mlabel{puresub} Suppose the Hecke algebra $\H^\G_\O$ is
  integrally quasi Frobenius. Let $\T:=\T^{\G}(\PR)$ be the 
$\H_{\O}^{\G}$-module 
$$
T=\bigoplus_{G'\in\L_{{G},L}}\sigma_{G'}\H_{\O}^{\G}.
$$ 
\begin{enumerate}
\item Let the functor $H_\O^\G=H_\O$ be defined as in
  \rmref{assochecke}. Then  
$$  
H_\O(\bX_h)=H_\O(\bX)=H_\O(\bX').
$$
\item Restricting endomorphisms to sublattices induces isomorphisms
$$
\End_{\O\G}(\bX_h)\cong\End_{\O\G}(\bX)\cong\End_{\O\G}(\bX').
$$
\item For every lattice $\bZ=\bX_h,\bX,\bX'$ the functor $H_\O$
  induces an isomorphism
$$  
H_\O:\End_{\H_\O^\G}(\T)\to\End_{\O\G}(\bZ).
$$
\end{enumerate}\end{thm}

\pf We may identify $\bX$ with its image $\bX_\brho$ under the
inejection $\brho$. We apply Theorem \ref{sublattices}. For any
$\O\G$-lattice $V$ we have always $H_\O(V)=\H_\O(\tau_{Q_L}(V))$ by
\cite[2.17]{di}, hence i) follows.  

By Lemma 
\ref{tauundtfunct} the endomorphisms of $\bX'$ map the subspaces 
$\tau_{Q_L}\bX'=\bX_h$ and $\tau_\bY(\bX')=\bX=\bX_\brho$ into itself,
restriction defines so an homomorphism from $\End_{\O\G}{\bX'}$ to 
$\End_{\O\G}(\bX_h)$ and to $\End_{\O\G}(\bX)$ repectively. Similarly
we have an homomorphism from $\End_{\O\G}(\bX)$ to
$\End_{\O\G}(\bX_h)$. By Theorem
\ref{hom1corr} part ii) restricting endomorphisms from $\bX'$ to $\bX_h$ is an
isomorphism between the corresponding endomorphism rings, therefore 
from $\bX'$ to $\bX$ as well. This shows part ii). Part iii) follows
now from the cooresponding part of \ref{hom1corr}.
\epf

In the following we call the $\H_{\O}^{\G}$-lattice $\T$ the
\emph{parabolic tensor space} of $\H_{\O}^{\G}$. As usual 
$k\otimes _{\O}\T=\bar{\T}$. The endomorphism
ring of $\T$ is denoted by $\S_{\O}=\S_{\O}^{\G}(\PR)$ and is called
the \emph{parabolic $q$-Schur algebra} of $\H_{\O}^{\G}$. The reason for
these names is the following: We shall apply these results to finite reductive
groups. Here $\H_{\O}^{\G}$ is always a Hecke algebra associated with
some reflection group $W$ extended possibly by an abelian group. The one
dimensional representation $\O\sigma _{G'}$ of $\H_{\O}^{{G}'}$ is a
representation analogous to the alternating representation of symmetric
groups. The space $\T$ is then isomorphic to the sum of modules
analogous to the permutation
representation of $W$ on its parabolic subgroups (respectively to the sum over
its parabolic subgroups corresponding to the Levi subgroups
$G'\in\L_{{G},L}$), which is again for symmetric groups $W$ a space which
is closely connected to tensor space. 
In general, endomorphism rings of such tensor spaces are now called $q$-Schur
algebras. In type $B$ for instance $q$-Schur algebras have been
defined based not only on parabolic subgroups but on the larger set of
reflection subgroups, (\cite{DJMa1} and \cite{DuS2}). With this notation 
we have now the following main result for the projective restriction system 
$\PR$ using Theorems \ref{decmatrix} and \ref{puresub}:

\begin{cor}\label{main1} Let $\PR$ be a projective restriction system 
such that the Hecke alhgebra $\H^\G_\O$ is integrally Frobenius. Then 
the decomposition matrix of the associated parabolic $q$-Schur algebra 
$\S_\O^\G(\PR)$ is part of the $\ell$-modular decomposition matrix of $G$.
\end{cor}

\begin{rem}\label{redundant} We remark that in general the Hecke
  algebras $\H^{G'}_R(\PR)$ and $\H^{G''}_R(\PR)$
for $G',G''\in \L_{{G},L}$ can be the same. This happens since
$\H^{\G}_R(\PR)$ is not a Hecke algebra defined over the Weyl
group $\tilde{W}$ of $\G$ but over some subquotient, namely
over the group $W(L,X_L)$ which is the stabilizer of $X_L$ in
${\cal N}_{\tilde{N}}(L)/L\cong {\cal N}_{\tilde{W}}(W_L)/W_L$, where 
$W_L:=(\tilde{N}\cap
L)/T$ is the Weyl group of $L$. Note that in this case the 
$\H_{\O}^{\G}$-modules
$\sigma_{G'}\H_{\O}^{\G}$ and $\sigma_{G''}\H_{\O}^{\G}$ are
isomorphic. Thus the resulting $q$-Schur algebras are Morita
equivalent, and we could remove one of the corresponding summands in
our summation 
$$
T_\O= \bigoplus_{M\in\L_{{G},L}}\sigma_M\H_{\O}^{\G}.
$$
\end{rem}

We have seen that part of the $\ell$-modular decomposition matrix of
$\G$ can be calculated by computing decomposition numbers of parabolic
$q$-Schur algebras, which in turn are endomorphism rings of $q$-tensor
space of Hecke algebras, provided the latter are integrally Frobenius.
In the next section we shall construct projective restriction systems
which satisfy this condition.

\section{Finite groups of Lie-Type}\mlabel{lie}
\begin{summ}\mlabel{notlie}{\rm
Our main applications of the theory developed in the
preceding sections are on finite groups of
Lie type. For any such group, say $G$, bold face letter
denotes its underlying algebraic group, that is, for example ${\bf G}$ will
denote the underlying algebraic group for which $G$ is the set of fixed points
of some Frobenius morphism (the latter will always be assumed to be known). 
There exist finite groups of Lie type which 
can be defined via various non-isomorphic algebraic groups. However, which
one is considered will always follow from the context. For finite groups of
Lie type, the theory of Deligne and Lusztig plays an important role. There,
the dual groups of finite groups of Lie type and the dual groups of 
connected algebraic groups are considered. We denote the dual group
of some finite group of Lie type or some connected algebraic group with the
same letter, added upper indices '$^*$'. For $L\in \L_G$ we view 
$L^*$ as Levi subgroup of $G^*$ in the usual way, after choosing an
isomorphism between the root data of $G$ and $G^*$ (see
\cite[4.2]{ca1}). We assume in the following that we have chosen
once and for all such an isomorphism for the actually considered group 
$G$. Then the identification of $L^*$ with a Levi subgroup of $G^*$
becomes canonical.

The algebraic groups we 
consider are defined over algebraic closures of finite fields. Therefore
let $\F_q$ denote the field of $q$ elements for some power $q$ of $p$ and let
$\bar{\F}_q$ be its algebraic closure. 

Finally, if 
$s$ and $z$ are elements of $G$, we write $s\sim_{G}z$ if $s$ is conjugate
to $z$ in $G$.

Let $ {\bG}$ be a connected reductive algebraic group over 
$\bar{\F}_q$. Let $F$ be a Frobenius automorphism of ${ {\bf G}}$ and
let $ {G}$ be the set of $F$-fixed points.
For simplicity of notation we assume that the
center of ${ {\bf G}}$ is connected, although most of the following 
could be formulated more generally. By general theory the center of
every Levi subgroup of ${ {\bf G}}$ is connected too. }
\end{summ}

First we want to summerize some well known facts.

\begin{summ}\mlabel{not2lie}{\rm
By the theory of Deligne and Lusztig,
the ordinary characters of $ {G}$ are distributed in pairwise disjoint 
series
${\cal E}( {G},z)$, called (rational) Lusztig series,
where $z$ runs through a set of representatives of the
conjugacy classes of semi simple elements of $ {G}^*$. 
For details see e.g. \cite[chapter 7]{ca1}.
The Lusztig series are defined by the Deligne-Lusztig operator, which
maps characters of generalized Levi subgroups of $G$ to
generalized characters of $G$. For our purpose 
it is enough to know that for $L\in \L_G$ this operator and
HC-induction coincide (\cite[7.4.4]{ca1}). 
In particular, the Lusztig series are unions of HC-series
and $\R _L^G{\cal E}(L,z)\subset {\cal E}( {G},z)$ for semisimple 
$z\in L^*$. We have therefore, using Frobenius reciprocity:

\begin{lemma}\label{luhc}
Two characters in different Lusztig series are contained in different 
HC-series. Moreover, if $\chi$ is contained in ${\cal E}( {G},z)$ and 
$\psi$ is a constituent of $\T^G_L\chi$, then 
$\psi$ is in ${\cal E}(L,(\tilde{z}))$ for some semisimple $\tilde{z}$ 
which is conjugate in $G^*$ to $z$. 
\end{lemma}

There exists a canonical isomorphism between the Weyl groups $W$ of
$G$ and $W^*$ of $G^*$ (by virtue of the fixed isomorphism between
their root data, see e.g. \cite[4.2.3]{ca1}). 
As the torus $T$ of $G$ is contained in every  
$M\in \L_G$, the Weyl group $W$ acts by conjugation on the set of Levi
subgroups of $G$ and on the set of pairs $(M,\chi)$, 
where $M\in \L_G$ and $\chi$
an irreducible character of $M$. We denote the stabilizer of $M$ under this 
action by $\stab_W(M)$. On the other hand
$W^*$, and thus $W$ (via the canonical isomorphism), acts on the set of Levi
subgroups of $G^*$ and on the set of pairs
$(M^*,(z))$, where $z\in M^*$ is semisimple and $(z)$ denotes its
conjugacy class in $M^*$. It follows immediately that $\stab_W(M)$ is
the stabilizer of $M^*$ as well; (the identification of $M^*$ with a
Levi subgroup of $G^*$ is just defined that way). Now let 
$z\in M^*$. Then by the definition of the Deligne-Lusztig operator 
(see \cite[7.2]{ca1}) we have
\begin{equation}\label{lusztigconjugation}
\E(M,z)^w=\E(M^{w},z^w).
\end{equation}
}
\end{summ}
\begin{summ}\mlabel{blockseries}{\rm 
It was shown by Brou\'e and Michel in
\cite[2.2]{brmi} that for an $\ell$-regular semisimple element 
$z\in  {G}^*$, the set
\begin{equation}\label{sbloecke}
\cE_{\ell}( {G},z)=\bigcup\limits_{t\in C_{ {G}^*}(z)_{\ell}}
\cE( {G},zt)
\end{equation}
forms a union of $\ell$-modular blocks, where $C_{ {G}^*}(z)_{\ell}$ denotes
the $\ell$-elements of $C_{ {G}^*}(z)$. Moreover, Geck and Hiss have 
shown in \cite{gehi} that the characters in $\cE( {G},z)$ form a 
basic
set of the blocks $\cE_{\ell}( {G},z)$ if $\ell$ is a good prime
for $ {G}$ (e. g. $\ell$ arbitrary if $ {G}$ is of type $A$ and 
$\ell$ odd if $ {G}$ is of type $B$, $C$ or $D$). That is
their reduction modulo
$\ell$ forms a $\Z$-basis of the additive group of generalized Brauer
characters in the blocks $\cE_{\ell}( {G},z)$. Therefore, the 
decomposition
matrix of the characters in $\cE( {G},z)$ determine uniquely the 
decomposition matrix of the whole blocks $\cE_{\ell}( {G},z)$. 
Moreover,
the latter can be derived from the former by decomposition of 
Deligne-Lusztig induced characters, as shown in \cite{gehi}.
An important property of $ {G}$ is that the Hecke algebra 
$\End_{K {G}}(\R_L^{ {G}}X)$ for some cuspidal irreducible $KL$-module
$X$ is always untwisted by \cite[4.23]{lubu}.}
\end{summ}

\begin{summ}\mlabel{gelfand}{\rm 
We want to exhibit a projective restriction system 
$\PR(X_G,Y_L)$ for $(L,\chi_L)$ and some $L\in\L_{G}$. An important
role for the
representations of $G$ plays a certain character, the
\emph{Gelfand-Graev} character $\Gamma=\Gamma_G$ of $G$. For its definition
and basic properties we refer to the standard literature,
e.g. \cite{ca1} and we list only a few facts on $\Gamma$ which will be
needed later. One of its main features is that many cuspidal
irreducible characters are constituents of $\Gamma$ and that all
irreducible constituents of $\Gamma$ occur with multiplicity one in
it. These constituents are called \emph{regular} characters. Each
Lusztig series contains precisely one regular character. Another
important property of $\Gamma$ is that its HC-restriction
$\T_{G'}^G(\Gamma_G)$ to the Levi subgroup $G'$ of $G$ gives the
Gelfand-Graev character $\Gamma_{G'}$ of $G'$. Moreover $\Gamma$
is induced from a linear character of the unipotent radical of a Borel
subgroup which is in particular $\ell$-regular. As a consequence there
is a unique projective $\O G$-lattice $P_G$ affording
$\Gamma$. Moreover,
\begin{equation}\label{hcrstgelfand}
\T_{G'}^GP_G=P_{G'}.
\end{equation}
Let $\chi$ be an irreducible constituent of $\Gamma_G$. Then using the
notation introduced in \ref{notquot} we conclude that $\pi^G_{\chi}(P_G)$ 
is a lattice $V$ affording the character $\chi$. Now
$$
\Hom_{KG}(P_G,\R_L^GKX_L) =\Hom_{KL}(\T_L^G\Gamma_G,KX_L)
=\Hom_{KL}(\Gamma_L,KX_L),
$$
where $X_L$ is a lattice affording $\chi_L$. Thus the multiplicity of
$\chi$ in $\R_L^G\chi_L$ is one if $\chi\in \E(G,z)$ and $\chi_L$ is 
regular, and zero otherwise.}
\end{summ}

\begin{summ}\mlabel{not2gelfand}{\rm 
For the following result we also need the fact, due to Hiss,
(\cite[4.6.1]{hi3}), that for a given block $\cB$ of $G$ the summand
of $P_G$ belonging to $\cB$ is indecomposable. In fact every union
\ref{sbloecke} of blocks of $G$ contains exactly one such
indecomposable direct summand of $P_G$. Let $M\in\L_G$ and let $\cB$ be
a union \ref{sbloecke} of blocks of $M$ containing the Lusztig series
$\E(M,y)$, where $y\in M^*$ is semi simple. Then we denote the unique
indecomposable direct summand of $P_M$ in $\cB$ by $Y_{M,(y)}$. Since
every Lusztig series contains a constituent of $P_M$, the summand
$Y_{M,(y)}$ is well defined, however the label $y$ is only 
determined up to its $\ell$-part by 
the indecomposable direct summand of $P_M$. In particular we have a
bijection between direct indecomposable summands of $P_M$ and
$M^*$-conjugacy classes of semi simple elements of order prime to
$\ell$.

Now consider the character $\psi$ of $\T^G_{G'}Y_{G,(y)}$ for some 
semisimple $\ell$-regular element $y\in G^*$. Then \ref{luhc} implies 
that $\psi$
consists precisely of those summands of $\Gamma_{G'}=\T^G_{G'}\Gamma_G$ 
which are in some
$\E(G',x)$ where $y$ is conjugate to $x$ in $G^*$ (in particular
$\T^G_{G'}Y_{M,(y)}=(0)$ in case that $G^{\prime *}$ 
contains no conjugate of $y$).
Now we can refine Equation \ref{hcrstgelfand} to  
\begin{equation}\mlabel{hcrstsumgelfand}
\T^G_{G'}Y_{G,(y)}=\bigoplus\limits_{x\in C_{G^{\prime *},(y)}}Y_{G',(x)},
\end{equation}
where $C_{G^{\prime *},y}$ is a set of representatives of the 
$G^{\prime *}$-conjugacy
classes making up the intersection of $G^{\prime *}$ with the 
$G^*$-conjugacy class of $y$. In particular, the indecomposable direct
summands of $\T^G_{G'}Y_{G,(y)} $ are contained in pairwise different blocks.}
\end{summ}

\begin{notation}
\mlabel{liegroupsnot}{\rm
In the following let $G$ be as above. Let $L\in {\cal L}_{{G}}$.
The underlying algebraic groups ${\bf G}$ and 
${\bf L}$ have connected centers.
We fix a semisimple element $z\in L^*$ and assume that we have a cuspidal
element $\chi_L$ in $\cE(L,z)$. 

For $M\in {\cal L}_{ G,L}$ such that $M^*$
  contains the semisimple Element $y\in G^*$ we set
  $X_{M,(y)}=\pi_{(L,\chi_L)}^M(Y_{M,(y)})$. 

Recall the action of $W$ on the set of pairs
$(M^*,(y))$, where $M^*\in \L_G$, $y\in M^*$ is semi simple and 
$(y)$ denotes its conjugacy class in $M^*$. Then $\stab_W(L^*)$ acts on
the set of $L^*$-conjugacy classes of semi simple elements.
In the following $\stab_W(z)$ denotes the stabilizer of $(z)$ under
this action. We denote the $\ell'$-part of any 
semisimple element $y\in G^*$ by $y'$.}\end{notation}

\begin{thm}\mlabel{centrellprime} Suppose that $\chi_L$ is a regular
character,  and that $\stab_W(z)=\stab_W(z')$. Then
$\PR(X_{G,(z)},Y_{L,(z)})$ is a projective restriction 
system for $(L,\chi_L)$ in $G$ consisting of the data 
$$
\PR(X_{G,(z)},Y_{L,(z)})=\{X_M,Y_M,X'_M\mid M\in\L_{G,L}\},
$$
where $X_M=X_{M,(z)}$, and $Y_M$ and $X'_M$ are defined as in  
\ref{defrestsys}.
\end{thm}

\pf By \ref{gelfand} we only have to show part (ii) of Definiton
\ref{defrestsys}.  

As usual we denote the character of $X_G=X_{G,(z)}$ by $\chi_G$.
Let $M\in {\cal L}_{G,L}$ and $\chi_1\not=\chi_2$ two summands 
of $\T^G_M\chi_G$.
By \ref{blockseries} we are done if we can show that $\chi_1$ and
$\chi_2$ are in two different unions of series defined in \ref{sbloecke} and 
hence in different blocks.
So assume that this 
is not the case. We may assume that both 
$\chi_1$ and $\chi_2$ are in the collection $\cE_{\ell}(M,z')$ 
of blocks of $M$.

By Lemma \ref{luhc}
all summands of $\T^G_L\chi_G$ and hence of $\T^M_L\chi_i$ ($i=1,2$)
are in series $\cE(L,y)$ for $y\in L^*$ conjugate to $z$ in
$G^*$. Moreover, by Equation \ref{lusztigconjugation} and general
HC-theory, both  $\T^M_L\chi_1$ and $\T^M_L\chi_2$ actually do have
summands in in $\cE_{\ell}(L,z')$ as the latter is the
union of the series $\cE(L,y)$ with $y'=z'$.
Thus $\T^G_L\chi_G$ has at least two summands in 
$\cE_{\ell}(L,z')$. However, by Mackey decomposition and cuspidality
of $\chi_L$, 
$$
\T^G_L\chi_G\le \R^G_L\T^G_L\chi_L=\bigoplus \limits_{x\in
  \stab_{W}(L)/(W\cap L)}\chi_L^x.
$$
By Frobenius reciprocity, $\chi_L$ has multiplicity one in
$\T^G_L\chi_G$, hence it follows that there exists $x\in \stab_W(L)$
such that 
\begin{equation}\mlabel{conjchil}
\chi_L\not=\chi_L^x\in  \cE_{\ell}(L,z').
\end{equation}
Now we extend Equation \ref{lusztigconjugation} to
the union of series $\cE_{\ell}(L,z')$. As the two unions of series 
$\cE_{\ell}(L,z')$ and $\cE_{\ell}(L,z^{\prime x})$ are either equal
or disjoint it follows 
$$
\cE_{\ell}(L,z')=\cE_{\ell}(L,z')^x=\cE_{\ell}(L,z^{\prime x}),
$$
hence $x\in \stab_W(z')$. By assumption we then also have $x\in \stab_W(z)$.
Thus 
$$
\chi_L^x\in \cE(L,z)^x=\cE(L,z^x)=\cE(L,z).
$$
However, by construction (see \cite[chapter 8]{ca1}), $\Gamma_L$ is invariant
under the action of $\stab_W(L)$, hence $\chi_L^x$ is regular. As
there is only one regular character in $\cE(L,z)$, namely $\chi_L$, we 
conclude $\chi_L^x=\chi_L$, a contradiction to \ref{conjchil}.
\epf 

Note that the assumption of the previuous theorem fits precisely the 
situation considered in \cite[section 5]{difl1}. There the Hecke-algebra part 
for $\R_L^GX_L$ was done, and our theorem here provides its extension to 
$q$-Schur algebras. For the general linear groups $G=GL_n(q)$ this was done in 
\cite{dija1}. In fact it was shown there, that we get a complete list of the 
irreducible $\ell$-modular representations applying our main results 
\ref{bhomfunctor},\ref{centrellprime} and Result \ref{1.3} to the Lusztig 
series $\cE(G,z)$, where $z$ runs through a set of representatives of 
$\ell$-regular conjugacy classes of $G=GL_n(q)$. Thus we obtain from 
Results \ref{lattiso} and \ref{lattiso2} a new proof 
for the following results from \cite{di2},\cite{di3},\cite{ja} and 
\cite{dija1} (compare \cite[4.15, 6.11]{di3}):    

\begin{thm}\label{oldstuff1} Let $G=GL_n(q)$. Then every cuspidal irreducible 
$KG$-module remains irreducible if reduced modulo $\ell$. On the other hand, 
every irreducible cuspidal $kG$-module is liftable to an irreducible 
$KG$-module. If $L$ is a Levi subgroup of $G$, and 
$X_L$ is an irreducible cuspidal $\O L$-lattice such that $\R_L^G(X_L)$ is 
reduction stable then the combined image $\tau_{P_G}(\R_L^GX_L)$ of the 
Gelfand-Graev representation in $\R_L^GX_L$ is pure in $\R_L^GX_L$. 
\end{thm}

These results can be generalized to classical groups under the additional 
assumption that the prime $\ell$ is linear. For details we refer to 
\cite{grhi}. Here is another situation, where our results apply:

\begin{cor}\label{appl2} Suppose that $\chi_L$ is regular and
$\bar{\chi}_L$ is irreducible. Moreover, suppose that the stabilizers
$\stab_W(\chi_L)$ and $\stab_W(\bar{\chi}_L)$ of $\chi_L$ and
$\bar{\chi}_L$ under the action of $\stab_W(L)$ are equal. Then 
$\PR(X_{G,(z)},Y_{L,(z)})$ is a projective restriction 
system for $(L,\chi_L)$.\end{cor}

\pf $\chi_L^x$ for $x\in \stab_W(L)$ is regular (see the last proof), 
hence it is the unique regular character in $\E(L,z^x)$ and 
$\stab_W(z)\le \stab_W(\chi_L)$. As the converse follows from Equation
\ref{lusztigconjugation}, we have 
equality. Since there is a bijection between the set of irreducible Brauer 
characters of $kL$ and representatives of
isomorphism classes of indecomposable projective $\O L$-lattices, we
see that $\stab_W(\bar{\chi}_L)$ is equal to the stabilizer of
$Y_{L,(z')}$ under the action of $\stab_W(L)$ (analogeously denoted by 
$\stab_W(Y_{L,(z')})$). As $P_L$ is invariant
under this action and has exactly one indecomposable direct summand in 
every union of series $\cE_{\ell}(L,z')$, we get
$\stab_W(Y_{L,(z')})=\stab_W(z')$ and the assertion follows from Theorem 
\ref{centrellprime}. 
\epf

In the situation of Theorem \ref{centrellprime} as well as of Corollary 
\ref{appl2} we have exhibited projective restriction systems $\PR(X_G,Y_L)$, 
which allow us to apply our main results of the previous section. In both 
cases the Hecke algebra $H=\End_{\O G}(\R_LX_L)$ is a Hecke algebra 
associated with an extension of a reflection group and it is known that it 
admits an associative (in fact symmetric)
bilinear form, whose determinant is a unit, hence is integrally Frobenius,  
(\cite{difl1} in the first and 
\cite{gehima2} in the second case). Thus Theorem \ref{main1} implies:

\begin{cor}\label{main2} The decomposition matrix of the $q$-Schur algebra 
$\S^G_\O=\S^G_\O(\PR)$ is part of the $\ell$-modular decomposition matrix 
of $G$.
\end{cor}

We assume now that $G$ is a Levi
subgroup of some larger group of Lie type $\G$. Again we assume that
the center of $\bf \G$ is connected. 

\begin{thm}\mlabel{indcuspgen}
Suppose that $\cE(L,z)$ has only one cuspidal element $\chi_L$
and that $z=z'$.
Let $\chi_G$ be a summand of $\R_L^G\chi_L$ of multiplicity one.
Suppose further that there exists some 
indecomposable projective $\O G$-lattice $Y_G$ such that
the character of $X_G=\pi_{(L,\chi_L)}^G(Y_G)$ is equal to $\chi_G$.
Then there exists an $\O L$-lattice $X_L$ 
with character $\chi_L$ and projective cover $Y_L$
such that $\PR(X_G,Y_L)$ is a projective restriction 
system for $(L,\chi_L)$ in $\G$.\end{thm}

\pf First we show part (i) of Definiton
\ref{defrestsys}. To do so we show that if $L^g=L$ with $g\in \G$
but $L^g\not\sim_G L$
then the character of $Y_G$ has no summand in common
with $\R_{L^g}^G\chi_L^g$. Recall that we can restrict attention to
standard Levi subgroups, hence we can assume $g\in \tilde{N}$,
respectively we can take $g=w$ in the Weyl group $W$ of $\G$. Then
by \ref{lusztigconjugation} 
$$
\chi_L^w\in \cE(L,z)^w=\cE(L^w,z^w)
$$
(using the action of $W$ on the conjugacy classes of $\G^*$ via
$W\simeq W^*$) and hence  
$$
\R_{L^g}^G\chi_L^g\in \cE(G,z^w)\subset \cE_{\ell}(G,z^w).
$$
However, $z=z'$, hence if $z^w$ is not conjugate to $z$ in $G$ then
we have for the set of blocks
$\cE_{\ell}(G,z^w)\cap \cE_{\ell}(G,z)=\emptyset$. As $Y_G$ is
indecomposable, we are done.

For part (ii) of Definiton \ref{defrestsys} we show that the summands 
of $\T^G_M\chi_G$ for $M\in \L_{G,L}$ are in pairwise different blocks.
So assume this is not the case, that is let $\chi_1,\chi_2$ be
summands of $\T^G_M\chi_G$ in the same block. Again we may assume
that one of them is $\chi_M\in \E(M,z)$. By \ref{not2lie}, every 
summand of $\T^G_M\chi_G$ lies in some series $\cE(M,y)\subset
\cE_{\ell}(M,y)$ for $y\in M^*$ conjugate to $z$ in $G^*$. 
As $z=z'$, it follows then by \cite[2.2]{brmi} that both $\chi_1$ and
$\chi_2$ are in $\E(M,z)$.
Like in the proof of Theorem \ref{centrellprime} we conclude that both 
of $\T^M_L\chi_i$ has a summand in $\E(L,z)$. However, by general
HC-theory, these summands must be cuspidal and thus by our assumption
they must be equal to $\chi_L$. Thus $\chi_L$ has multiplicity of at
least two in $\T^G_L\chi_G$. By Frobenius reciprocity this gives a
contradiction to our assumption. 
\epf

\begin{rem}\mlabel{onechar}\rm If $L$ is a group of classical type, every 
Lusztig
series has at most one cuspidal character (see \cite[13]{ca1}).\end{rem}

The following corollary now follows immediately from the preceeding sections
\begin{cor}\mlabel{indcuspgencor}
Let the assumptions be as in Theorem \rmref{indcuspgen}.
Then the $q$-Schur algebra 
${\S}_{\O}^{ {\G}}(\PR(X_G,Y_L))$ 
is defined and its decomposition matrix is a submatrix of the 
decomposition matrix of $\cE(\G,z)$.\end{cor}

An example of particular interest for the previous two Theorems 
are the regular characters in series $\cE(L,y)$ for $y=y'$. We call
such series $\ell$-regular Lusztig series. Note that if some 
Lusztig series 
has just one element then this element is a regular character. 


There are two main examples of the above, which have been investigated 
in previous papers by various authors. We keep the outline of the
following examples sketchy, for further detailes the reader may refer
to the cited references.
\begin{summ}\mlabel{lin}{\rm 
The first one is
the well known theory for finite general linear groups. It was developed
in several articles by G. James and the first named author
(see \cite{di}, \cite{di2},
\cite{di3}, \cite{dija1}, \cite{ja} and \cite{ja1}). First one observes
that all the Hecke algebras appearing
as endomorphism rings of induced irreducible cuspidal $KL$-modules for
$L\in {\cal L}_{{G}}$,
are defined over groups $H$ isomorphic to direct products of symmetric groups.
Secondly, any irreducible cuspidal $KL$-module for $L\in \L_G$ is regular.
Therefore one can apply Theorem \ref{centrellprime} to cuspidal
characters in $\ell$-regular Lusztig series. It follows that the 
decomposition matrix of the $q$-Schur algebra
${\cal S}_{\O}(\PR(X_{{G},(z)},Y_{L,(z)}))$ gives a part of the
decomposition matrix of $\cE(G,z)$, where $z$ is $\ell$-regular and $L\in
\L_G$ is a minimal Levi subgroup containing $z$ (for $G$ a general
linear group, we have $G\cong G^*$ canonically). 
Moreover, it turns out
that ${\cal S}_{\O}(\PR(X_{{G},(z)},Y_{L,(z)}))$ 
has the same number of isomorphism classes of projective indecomposable
modules as ${\cal S}_{K}(\PR(X_{{G},(z)},X_{L,(z)}))$ 
(see \cite{dija1}). 
However, the latter are by definition
in bijection to the elements $S({G}/L,\chi_L)$, where $\chi_L$ is the unique
(cuspidal) character in $\cE(L,z)$. As HC-series and 
Lusztig series for finite general linear groups coincide,
it follows that the decomposition numbers of the various 
$q$-Schur algebras defined over cuspidal  irreducible regular characters in
$\ell$-regular Lusztig series gives us the complete decomposition matrix
of all series $\cE({G},z)$ for $\ell$-regular elements $z\in G^*=G$.
The decomposition matrix of the other characters can now easily deduced
by decomposition of Deligne Lusztig induced characters. However,
the results in this article can be applied to an arbitrary group ${G}$ 
such that
$\bf G$ is of type $A$ with connected center. 
Thus it follows that the decomposition numbers of 
the $\ell$-regular Lusztig series of ${G}$ are the same as the
decomposition numbers of $\ell$-regular Lusztig series of suitable 
general linear groups. }\end{summ}

\begin{summ}\mlabel{linclass}{\rm 
The second example concernes other 
finite classical groups $ {\G}$ 
such as symplectic, special orthogonal and general unitary groups.
They were investigated in \cite{grhi} by G. Hiss and the second named
author.  
If one restricts attention to so called linear primes for these groups
(see \cite{grhi}, Section 1 for definition, and note that linear primes are 
roughly half of all primes dividing the order of $ {\G}$), 
one can develop a theory similar to the one of finite general linear groups.

First one observes that the reduction modulo $\ell$ of
every cuspidal character $\chi_L$ of a 
Levi subgroup $L$ of $ {\G}$ is irreducible. 

Secondly, let $\phi\in \cE({L},z), \psi \in \cE({M},s)$ for $L,M\in
\L_{ {\G}}$ and $z,s\in  {\G}^*$ $\ell$-regular, be two
cuspidal irreducible characters not conjugate in $ {N}$. Then 
two characters $\alpha\in S( {\G}/L,\phi)$ and $\beta\in
S( {\G}/M,\psi)$ are lying in different blocks. 

Thirdly, for every
cuspidal irreducible character in an $\ell$-regular Lusztig series ${\cal
  E}(L,z)$ for some $L\in \L_{ {\G}}$ one can
define a projective restriction system and an associated $q$-Schur algebra of
type $B$ or $D$.
Using the representation
theory of Hecke algebras over classical extended Weyl groups for
linear primes (see \cite{dija} and \cite{grhi}, Section 7) one again
finds that the $q$-Schur algebra ${\cal S}_\O (\PR(X_G,Y_L))$
has the same number of isomorphism classes of projective indecomposable
modules as the $q$-Schur algebra over $K$ (see \cite{grhi} Corollary 8.2).
Thus we again have the complete decomposition matrix
of all series $\cE( {\G},z)$ for $\ell$-regular elements $z\in
\G^*$.

In these examples, $Y_L$ is the unique indecomposable projective 
$\O L$-lattice affording
a character with some irreducible cuspidal constituent 
$\chi_L\in\cE(L,z)$,
while $X_G$ is a suitably chosen $\O G$-lattice affording a character $\chi_G$
with multiplicity one in 
$\R^G_L\chi_L$. $G\in \L_{ {\G}}$ is chosen as follows:
The Dynkin diagram of the Weyl group of $L$ has at most one connected
component $I$ of type $B$ or $D$ (if it has no such component we set
$I=\emptyset$). Then $G$ is a maximal Levi subgroup of $ {\G}$ containing
$L$ such that its Weyl group has one connected component of type $A$ and one
equal to $I$.}\end{summ}

We note that in both examples one can apply Corollary \ref{lattiso} to the
projective restriction systems associated to irreducible cuspidal characters
in $\ell$-regular Lusztig series. We get a complete description of the module
structure of the so called Steinberg lattice, that is the unique quotient
lattice of the Gelfand-Graev lattice affording the Steinberg character. In case
of general linear groups, the corresponding quotient of the Gelfand-Graev 
lattice
affording other irrducible characters (even for those characters not lying in
$\ell$-regular Lusztig series) are described in \cite{gr2}. There, also for
the special linear groups corresponding results have been obtained. 

\begin{summ}\mlabel{genunipot}{\rm 
Let $ {\bf \G}$ be an arbitrary reductive algebraic group
with connected center. For simplicity of notation we assume that the 
Dynkin diagram of $ {\bf \G}$ is connected. However, the following can 
be formulated more generally. Let $L\in {\cal L}_{ {\G}}$
have a unique cuspidal character $\chi_L$ in the series 
$\cE(L,z)$ for some $z\in L^*$. Let $\hat{X}$ be a $KL$-module with 
character $\chi_L$.
The Dynkin diagram of $L$ has at most one connected component ${\cal I}$
which is not of type $A$. Let ${\cal J}$ be the union of the connected
components of the Dynkin diagram of $L$ of type $A$. Then there exists
a maximal Levi subgroup
$G\in {\cal L}_{ {\G}}$ whose Dynkin diagram has at most two connected 
components, one being of type $A$ and having ${\cal J}$ as subset and one 
equal to ${\cal I}$. Then $\stab_W(\chi_L)$ is a direct product of symmetric
groups. We take $\chi_G$ to be the character in $S(G/L,\chi_L)$ corresponding
to the alternating character of $\stab_W(\chi_L)$. In particular,
$\pi_{(L,\chi_L)}(\T^G_M\chi_G)=\chi_M$ is irreducible for $
M\in {\cal L}_{G,L}$ and 
$\pi_{(L,\chi_L)}(\T^G_L\chi_G)=\chi_L$ by \cite{hole}, Theorem 5.9.
Finally suppose that there exists an indecomposable projective $\O G$-lattice
$Y_G$ such that $\pi_{(L,\chi_L)}(Y_G)=X_G$ has character $\chi_G$.
Setting $X_M=\pi_{(L,\chi_L)}(\T^G_MX_G)$ for $M\in {\cal L}_{G,L}$ and
denoting the projctive cover of $X_M$ by $Y_M$, we
can apply Theorem \ref{indcuspgen} to get a projective restriction
system $\PR(X_G,Y_L)$. Moreover, because of our assumptions,
we can apply Corollary \ref{indcuspgencor}. Obviously, if 
$L\cong G'\times L_1$ and $G\cong G'\times L_2$, where $G'$ is a 
finite group of Lie type with Dynkin diagram ${\cal I}$ (as it holds 
in the case of finite classical groups), it follows from
the theory of linear groups that $Y_G$ exists. Moreover, if $L$ is of 
classical type, there always exists at most one cuspidal character in a 
Lusztig series of $L$. We can formulate the following}\end{summ}

\begin{cor}\mlabel{classgeneral} 
Let the notation and assumptions be as in \ref{genunipot}.
Then the $q$-Schur algebra ${\S}_{\O}(\PR(X_G,Y_L))$ is
defined and its decomposition matrix embeds into the decomposition matrix of
$\cE( {\G},z)$.\end{cor}

We want to discuss an application of our theory that has not been investigated
in other articles, yet, the decomposition matrix of the
unipotent characters of the general unitary groups for arbitrary odd primes. 
To do so we need some
preliminary considerations about the general linear groups.

\begin{summ}\mlabel{linlabel}{\rm
For the following Lemma the reader might recall the labeling of 
irreducible unipotent characters of the general linear group $G=GL_n(q^2)$ by
partitions $\alpha \vdash n$ of $n$
(see \cite[13.8]{ca1}). Here, a unipotent Brauer
character is a Brauer character appearing as summand of the reduction 
modulo $\ell$ of some unipotent character. Using the lexicographic 
order for partitions (see \cite{ja}), the Steinberg character has the lowest
label. For any prime $\ell\not=p$, the decomposition matrix of 
the unipotent characters, arranged with respect to the lexicographic order
along a vertical edge of the decomposition matrix,
can be made unitriangular if the unipotent Brauer characters are 
arranged suitably on a horizontal edge of the decomposition matrix
(see again \cite{ja}). Moreover, the decomposition matrix of 
the unipotent characters is a square matrix. Giving the $i^{th}$ unipotent
Brauer character the label of the $i^{th}$ unipotent character, we
get a labeling of the unipotent Brauer characters by the partitions of $n$.
We remark that square unitriangularity of the decomposition matrix determines
the labeling of the Brauer characters uniquely. This follows easily from
the fact that a different labeling of the Brauer characters has the
effect of multiplying the decomposition matrix with a permutation matrix
and if the resulting matrix is again square unitriangular then the 
permutation matrix is necessarily the unit matrix.} \end{summ}

\begin{lemma}\mlabel{decnumlin}  
Let $\mu\not=(1^n) $ be a partition of $n$. Then there exists an odd prime
$\ell$ and some power $q$ of $p\not=\ell$ such that the 
unipotent Brauer character $\phi_{\mu}$ of $GL_n(q^2)$ with label $\mu$
appears in the reduction modulo
$\ell$ of at least two different irreducible unipotent characters of
$G$. Moreover, $\ell$ can be chosen to be linear for $GU_n(q)$.\end{lemma}

\pf The proof proceeds in two steps.
First we prove the assertion for
$\mu=(j,1^{n-j})$ for $j\ge 1$. If $n$ is odd, take $\ell$ to be 
a Zsigmondy prime of $q^n-1$ (see \cite[8.3]{hup}) for some $q$ an
arbitrary power of some prime $p$. The order of 
$q^2$ modulo $\ell$ is the same as the order of $q$ modulo $\ell$ and 
$\ell$ is linear for $GU_n(q)$. Now the assertion follows in this 
case immediately from \cite[6.5]{ja} by taking $e=n$ in that theorem.
For $n=2$ the  the assertion follows from \cite[6.5]{ja}
using $\ell=3$ and $q=4$. Now let $n>2$ be even. Let $q$ be an
arbitrary power of some prime $p$ (in case $n=6$ take $q\not=2$). Let $\ell$
be a Zsigmondy prime of $q{^{2n}}-1$. Then $\ell$ is again linear for $GU_n(q)$
and the assertion follows in this case from \cite[6.5]{ja} by taking 
$e=2n$.\\
The second step is to proceed by induction. 
Let $\mu=(\mu_1,\ldots, \mu_r)$. We set
$\mu'=(\mu_2,\ldots,\mu_r)\vdash n-\mu_1$. We can assume that 
$(0)\not=\mu'\not= (1^{n-\mu_1})$. 
By induction hypothesis there exists $\ell\not=p$, a power $q$ of $p$ and a
partition $\lambda'\not=\mu'$ of $n-\mu_1$ such that the $\ell$-modular 
unipotent Brauer character 
of $GL_{n-\mu_1}(q^2)$ with label $\mu'$ is a constituent of the
reduction modulo $\ell$ of the unipotent $KGL_{n-\mu_1}(q^2)$-character
with label $\lambda'=(\lambda_2,\ldots ,\lambda_s)$ and $\ell$ is 
linear for $GU_{n-\mu_1}(q)$. Moreover, by unitriangularity of the 
decomposition matrix it follows that $\lambda'\le \mu'$ in the 
lexicographical order.
Thus $\mu_1\ge \lambda _i$ for $2\le i\le s$. By \cite[6.18]{ja}
$\phi_{\mu}$ is a constituent of the reduction modulo $\ell$ of the unipotent 
$KGL_{n}(q)$-character with label $(\mu_1,\lambda_2,\ldots,\lambda_s)$.\epf

\begin{summ}\mlabel{unitary}{\rm
For the remainder of this section let $G$ be the finite general 
unitary group over the field with $q^2$ elements
and let $\ell\not=2$. We will investigate unipotent characters,
thus we assume that $\chi_L$ is a cuspidal unipotent irreducible character 
for some $L\in {\cal L}_{G}$. We remark however that the following 
can also be applied to some other $\ell$-regular Lusztig series ${\cal
  E}(G,z)$ by 
\ref{genunipot}, as 
the Levi subgroups of $G$ are direct products of at most one unitary
factor and general linear groups over $\F_{q^2}$. The characters in ${\cal
  E}(G,z)$ are canonically labeled by the unipotent characters of $C_{G}(z)$
and their decomposition matrices are the same (see \cite{grhi}).
By \cite[13.8]{ca1} the irreducible unipotent $KG$-characters
can be labeled by partitions of $n$. Similarly to the linear groups, it was
shown in \cite[6.6]{gehima1} that the decomposition matrix of $G$
is square lower unitriangular (for appropriate arrangements of the characters).
We again get a labeling of the unipotent Brauer characters by partitions
of $n$ like in the case of $GL_n(q)$. 
Now the unipotent characters of $G$ (and thus the Brauer characters)
can also be labeled by the bipartitions of $\{(n-r)/2\}_r$, where $r$
runs through all positive integers less or equal to $n$ such that
$n-r$ is even and the unitary group of degree $r$ has a unipotent cuspidal
character (see \cite{grhi}, for the relevant $r$ see
\cite[13.7]{ca1}). 
Here a bipartition of $n$ is an ordered pair of partitions 
$\alpha \vdash a$ and $\beta\vdash b$ with $a+b=n$. 
A unipotent character
has as label some bipartition of $(n-r)/2$ if and only if it lies in the 
HC-series of the cuspidal irreducible unipotent character 
$\chi_L$ of the Levi subgroup $L\cong GU_r(q)\times GL_1(q^2)^{(n-2)/2}$ of
$G$ (see \cite[13.8]{ca1}). }\end{summ}

Now consider $S(G/L,\chi_L)$. Then we have 
\begin{thm}
\mlabel{orderun}  
Let $\psi$ be the character in $S(G/L,\chi_L)$ with least label in
the lexicographical order.
Then $\psi$ appears with multiplicity one in $\R_L^{G}\chi$.\end{thm}

\pf We turn to the labeling of the characters in 
$S(G/L,\chi)$ by the bipartitions of $s=(n-r)/2$. Suppose that
$\psi$ has label $(\alpha,\beta)$ for partitions $\alpha\vdash a$ and 
$\beta\vdash b=s-a$. By the assumption on $\psi$, it follows from 
the lower unitriangularity of the decomposition matrices 
that there exists no odd prime $\ell'$ and no power $q'$ of some prime $p'$ 
different to $\ell'$ such that the $\ell'$-modular irreducible unipotent 
Brauer character with the same label as $\psi$ appears as constituent of the 
reduction modulo $\ell$ of some element in $S(G'/L',\chi')$ different
to $\psi$. Here $G'$ denotes $GU_n(q')$, $L'$ the a 
standard Levi subgroup of $G'$ isomorphic to
$GU_r(q')\times GL_1({q'}^2)^{(n-2)/2}$ and $\chi'$ the cuspidal irreducible
unipotent character of $G'$. Recall that by general
theory the multiplicities of summands $\rho'$ in $\R_{L'}^{G'} \chi'$
and $\rho$ in $\R_L^G\chi_L$ are the same, if $\rho$ and $\rho'$ have
the same label. Now assume that
$\alpha\not=(1^a)$. Then choose $\ell'$ and $q'$ with 
Lemma \ref{decnumlin} such that the $\ell'$-modular
irreducible unipotent Brauer character of 
$GL_{a}({q'}^2)$ with label $\alpha$ appears in the reduction modulo $\ell'$
of some irreducible unipotent character with label $\gamma\not=\alpha$ and
$\ell'$ is linear for $GU_a(q')$. Linearity of $\ell'$ only depends on $q'$,
thus $\ell'$ is linear for $GU_n(q')$. Now it follows from \cite[3.36]{gr} 
that the $\ell$-modular unipotent Brauer character of $GU_n(q')$ is a 
constituent of the reduction modulo $\ell$ of the character in 
$S(G'/L',\chi_L')$ with label $(\gamma,\delta)$, a contradiction to
our statement at the beginning of the proof. Thus $ \alpha=(1^a)$. 
Similarly one shows that $\beta=(1^{b})$.

Now recall the way how one passes from labeling of the characters in
$\cE(G,1)$ by partitions of $n$ to their labeling by bipartitions.
This way is described for example in \cite{fosribl}. If $r=0$,
the partition $(1^n)$ is the label of the Steinberg character and the least
partition in the lexicographical order. Thus we may assume that $r>0$.
Assume that $a>0$ and $b>0$. Then one can see immediately, that if $\lambda$ 
is the partition labeling $\psi$, then $\lambda$ must be of one of the 
following forms.
\[a) \quad (m+3+b,m+3+(b-1),\ldots,m+3,m,m-1,\ldots,2,1^{2a+1}),\]
\[b) \quad (2+b,2+(b-1),\ldots,3,1^{2a}),\]
\[c) \quad (3+b,3+(b-1),\ldots,4,1^{2a+1}),\]
for $m\ge 2$. However, consider the partition 
$\delta=(m+b+1,\ldots,2,1^{2a+2b+1})$ in case $a)$,
respectively $\delta=(b,\ldots,2,1^{2a+2b+1})$ in cases $b)$ and $c)$. Then
$\delta$ is a partition of $n$ and corresponds to a bipartition of the form 
$((0),(1^s))$ or $((1^s),(0))$. Moreover, $\delta$ is lower in the 
lexicographical order than $\lambda$, a contradiction to our
assumption on $\psi$. Thus it follows that $a=0$ or $b=0$. Now the assertion 
follows from the representation theory of the Weyl group of type $B$, see
\cite{ca1}. \epf
 
As a consequence of Proposition \ref{orderun} we can apply Theorem
\ref{indcuspgen} to get our last theorem. We use the following notation
for it. Let $\chi_L$ be an irreducible cupidal unipotent character of 
$L\in {\cal L}_{G}$. Let $\chi_G$ be the character of lowest label in
$S(G/L,\chi_L)$. Let $\phi$ be the unipotent Brauer character of
$G$ with the same label as $\chi_G$ and let $Y_{G}$ be an indecomposable 
projective $\O G$-lattice such that $Y_G/\Jac{Y_G}$ has Brauer
character $\phi$.
Let $X_{M}=\pi_{(L,\chi_L)}(\T^G_MY_G)$ for $M\in {\cal L}_{G,L}$, let
$\chi_M$ be its character and let $Y_M$ be the projective cover of $X_{M}$. 
 
\begin{thm}\mlabel{schurunipotunit} The data 
$\{(\chi_M,X_M,Y_M)\mid M\in {\cal L}_{G,L}\}$ define a projective 
restriction system $\PR(X_G,Y_L)$. 
In particular, the decomposition matrix of
${\S}_{\O}(\PR(X_G,Y_L))$ is a submatrix of the decomposition matrix of
$\cE(G,1)$.\end{thm}

\pf By square unitriangularity of the 
decomposition matrix of $\cE(G,1)$
it follows that $\phi$ is not a constituent of the reduction modulo $\ell$
of any character in $S(G/L,\chi_L)$ different to $\chi_L$ and its
multiplicity in $\bar{\chi}_L$ is one. Thus 
$\pi_{(L,\chi_L)}(Y_G)=X_G$ has character $\chi_G$. Now the first statement
follows from Theorem \ref{indcuspgen}. The second
statement follows from Corollary \ref{indcuspgencor}.\epf

\begin{summ}\mlabel{sln}{\rm
We have so far only considered finite Lie groups whose underlying
algebraic group has
connected center. In case of classical groups it was shown in \cite{grhi} that
one can extend the results for linear primes to groups whose underlying
algebraic group has non-connected center. In general this is done by embedding
the group in a finite Lie groups whose underlying algebraic group has
conncted center. In case of classical groups the center of the underlying
algebraic group has two connected components. We get more problems in case the
number of connected components increase. As a standard example for this to
happen we can consider the finite special linear groups. These groups were
treated in \cite{gr2}. There it was shown how one can describe the
decomposition matrix of a set of irreducible characters whose reduction modulo
$\ell$ generate the group of generalized Brauer characters in terms of
decomposition matrices of $q$-Schur algebras defined over extended Weyl groups
of type $A$.}\end{summ}

\footnotesize\baselineskip6pt

\end{document}